\documentclass[11pt]{article}
\usepackage{graphicx}
\usepackage{color}
\usepackage{amsmath,amsfonts,amssymb}
\def\la{\big\langle}
\def\ra{\big\rangle}

\def\ds{\displaystyle}
\def\forall{\hbox{for all}~}
\def\L{{\bf L}}

\def\bfv{{\bf v}}

\def\bfw{{\bf w}}
\def\bfu{{\bf u}}

\def\bfn{{\bf n}}

\def\bff{{\bf f}}
\def\bfB{{\bf B}}
\def\bfH{{\bf H}}
\def\bfF{{\bf F}}
\def\bfG{{\bf G}}

\def\bfA{{\bf A}}

\def\bpm{\begin{pmatrix}}
\def\epm{\end{pmatrix}}

\def\ve{\varepsilon}

\def\E{{\cal E}}

\def\R{{\mathbb R}}

\def\vp{\varphi}

\def\vs{\vskip 2em}

\def\v{\vskip 1em}
\def\D{{\cal D}}
\def\O{{\cal O}}

\def\begi{\begin{itemize}}
\def\endi{\end{itemize}}
\def\C{{\cal C}}
\def\div{\hbox{div}}
\def\ov{\overline}
\def\Tilde{\widetilde}
\def\Hat{\widehat}
\def\bega{\begin{array}}
\def\enda{\end{array}}

\def\bel{\begin{equation}\label}
\def\eeq{\end{equation}}
\def\sqr#1#2{\vbox{\hrule height .#2pt
\hbox{\vrule width .#2pt height #1pt \kern #1pt
\vrule width .#2pt}\hrule height .#2pt }}
\def\square{\sqr74}
\def\endproof{\hphantom{MM}\hfill\llap{$\square$}\goodbreak}

\newtheorem{theorem}{Theorem}[section]

\newtheorem{lemma}{Lemma}[section]

\newtheorem{remark}{Remark}[section]

\begin{document}
\title{\bf  The Initial Stages of a Generic Singularity \\  for a 2D Pressureless Gas}
\vs

\author{Alberto Bressan$^{*}$, Geng Chen$^{**}$, and Shoujun Huang$^{\dag}$\\
\, \\
$^*$Department of Mathematics, Penn State University, \\
University Park, PA ~16802, USA.\\
$^{**}$Department of Mathematics, University of Kansas,\\
 Lawrence, KS 66045, USA.\\
$^\dag$College of Mathematical Medicine, Zhejiang Normal University, \\
Jinhua 321004, P. R. China.\\
\, \\
E-mails: axb62@psu.edu,~gengchen@ku.edu,~sjhuang@zjnu.edu.cn.}
\maketitle
\v
{\bf Abstract:}  We consider the Cauchy problem for the equations of pressureless gases in two space dimensions.   For a generic set of smooth initial data (density and velocity), it is known that
the solution loses regularity at a finite time $t_0$, where both the  the density and the velocity gradient become unbounded.
Aim of  this paper is to provide an asymptotic description of the solution beyond the time of
singularity formation.    For $t>t_0$ we show that a singular curve is formed, where the mass has positive density w.r.t.~1-dimensional Hausdorff measure.
The system of equations describing the behavior of the singular curve is not hyperbolic.
Working within a class of analytic data, local solutions can be constructed
using a version of the Cauchy-Kovalevskaya theorem. For this purpose, by a suitable change of variables we rewrite the evolution equations as a first order system of Briot-Bouquet type, to which a general existence-uniqueness theorem can then be applied.


\section{Introduction}
\label{s:1}
\setcounter{equation}{0}

We consider the initial value problem for the equations of pressureless gases
in two space dimensions:
\begin{equation}\label{SP}
\left\{\begin{array}{rl}
\!\!\partial_t\rho+\div( \rho \bfv)&=0, \cr
\!\!\partial_t(\rho \bfv)+\div\big( \rho \bfv\otimes \bfv\big)&=0,\enda
\right.  \qquad\quad  t\in  ]0,T[\,,\quad x\in\R^2,\eeq
\bel{ic}\rho(0,x)~=~\bar \rho(x), \qquad \bfv(0,x)~=~\bfw(x).\eeq

As long as the solution remains smooth, it can be directly computed by the method
of characteristics.  Indeed,
\bel{vel} \bfv(x+ t\bfw(x))~=~\bfw(x).\eeq
In addition, conservation of mass implies
\bel{den}\rho(x+t\bfw(x))~=~{\bar \rho(x)\over \det\bigl(I + tD\bfw(x)\bigr)}\,.\eeq
The equations (\ref{vel})-(\ref{den}) describe the solution up to the first time
$t_0$
when a singularity occurs.    This happens as soon as one of the matrices
$I+ tD\bfw(x)$ is no longer invertible.  Equivalently,
$t_0$ can be characterized as the smallest time $t>0$
such that $-1/t$ is an eigenvalue of $D\bfw(x)$, for some $x\in\R^2$.
%

In one space dimension, the solution will generically contain a finite number of point masses,
which can eventually collide with each other \cite{D}.
In two space dimensions, the equations (\ref{SP}) give rise to
singular curves, where the mass has positive density w.r.t.~1-dimensional
Hausdorff measure.   The time evolution of these singular curves has been
described in \cite{BCH}, together with their generic interaction patterns.
 An important aspect of
the governing equations for these curves is that they are not of hyperbolic type.
Therefore, apart from some special configurations, the Cauchy problem is ill posed.
Local solutions can  be constructed only within a class of analytic data,
as in the case of vortex sheets \cite{CO}.

In the present paper, our main goal is to provide a detailed asymptotic description of
how a new singular curve is formed.  We thus consider analytic initial data
$\bar\rho, \bfw$, defined on an open set $V\subset\R^2$, satisfying the following assumptions.

\begi
\item[{\bf (A1)}] {\it For every $x\in V$, the initial density is strictly positive: $\bar\rho(x)>0$.
Moreover, the Jacobian matrix $D\bfw$
of the  initial velocity  has real distinct eigenvectors: $\lambda_1(x)<\lambda_2(x)$.}

\item[{\bf (A2)}]  {\it There exists a point $x_0\in V$ where the smaller eigenvalue $\lambda_1(x_0)<0$
attains a negative minimum.
At this point,  the Hessian matrix
$D^2\lambda_1(x_0)$ is strictly positive defined.}
\endi
Calling
\bel{taudef} \tau(x)~\doteq~{-1\over \lambda_1(x)}\eeq
the time when the determinant vanishes: $\det(I+tD\bfw(x)\bigr)=0$,
by {\bf (A2))} it follows that $\tau$ has a positive minimum at $x_0$, and
the Hessian matrix $D^2 \tau(x_0)$ is also strictly positive definite.

The remainder of the paper is organized as follows. In
Section~\ref{s:2} we review the system of equations describing
a singular curve \cite{AR, BCH}.

As a ``warm up", Section~\ref{s:3} provides an asymptotic description
of the emergence of singularities in dimension 1. In some ways, this resembles the
formation of a new shock for a scalar conservation law.   The analysis relies on
tools from differential geometry originally developed in
 \cite{GS, Gk2}.   As shown in Figures~\ref{f:gen52} and \ref{f:gen51}, there exists a
 cusp region $\Omega$ in the $t$-$y$ plane whose points can be reached by three different
 characteristics, starting at points $x_1<x_2<x_3$.  For each $(t,y)\in \Omega$, these
 points are found by solving a cubic-like equation.   By parameterizing the two extremal
 sheets in terms of the middle sheet, so that $x_1=x_1(t, x_2)$ and $x_3 = x_3(t, x_2)$,
the problem is reduced to the study of a quadratic equation.  In the end, the singular curve
is described in terms of the stable manifold of a vector field $\bff$
with smooth coefficients in the $t$-$x$ plane, see (\ref{vf}).

The last two sections deal with
singularity formation in the 2-dimensional case.    In
Section~\ref{s:4} we begin by a detailed
analysis of the smooth solution provided  by the method of characteristics
at times $t>\tau(x_0)\doteq{-1\over \lambda_1(x_0)}$.  As shown in Fig.~\ref{f:gen42},
for $x$ in a neighborhood of $x_0$ the
 map $x\mapsto x+ t \bfw(x)$ is no longer one-to-one.   Indeed, as in the 1-dimensional case, points $y$ in a neighborhood
 of $x_0 + t \bfw(x_0)$ have three pre-images.  The two extreme ones are relevant for our
 construction.
 We start with a  leading order expansion and then construct
 the exact solutions by the Newton-Kantorovich iteration procedure.
 This yields a detailed asymptotic expansion of the velocity $\bfv^\pm$ and of the density $\rho^\pm$
 of the pressureless gas in front and behind the singular curve.
 Finally, in Section~\ref{s:5} we use the previous analysis to achieve a
construction of the singular curve during its initial stages, valid for analytic data.
Here the heart of the matter is to show that, by a suitable change of variables, the
evolution equations can be written as a first order system of Briot-Bouquet type, to which
a general existence-uniqueness theorem can be applied \cite{bao, li}.

The present results can be seen as part of a general research program aimed at understanding
the singularities of solutions to nonlinear PDEs, for generic data~\cite{Damon, Guck}.
In one space dimension,
singular solutions to the sticky particle model were studied in \cite{D}, while the recent paper
\cite{AR} derives the evolution equations for a singular surface in multidimensional space.
Generic singularities of solutions to a second order
variational wave equation were studied in \cite{BC, BHY, CCD}.
In the case of scalar conservation laws in one space dimension, properties of
generic solutions were studied
in \cite{GS, Gk2, S}.
In one space dimension, the analysis in \cite{Kong} provided  an asymptotic description
of this blow up, and also of the initial stages of the new shock, in a generic setting.
A study of shock formation for multidimensional scalar conservation laws can be found in
\cite{DM}.
For the 3-dimensional equations of isentropic gas dynamics,
the asymptotic of blow up, leading to shock formation, was studied in \cite{BSV}.
For earlier results on blowup for hyperbolic equations we refer to~\cite{A}.

The equations of pressureless gases have found applications as models
for the motion of granular media \cite{ERS}.
In particular, they provide
an approximate description of the large-scale distribution of matter in the universe
\cite{GSS, Z}.
More recently,  they have been shown to describe a limit of collective dynamics with short-range interactions \cite{ST}.

\section{Evolution of a singular curve}
\label{s:2}
\setcounter{equation}{0}

Following \cite{BCH}, we write a system of equations describing the evolution of a singular curve $\gamma$ in 2-dimensional space, as shown in
Fig.~\ref{f:gen77}.
At a given time $t$,
 the curve will be parameterized as
$$ \zeta~\mapsto ~y(t, \zeta),\qquad\qquad  \zeta\in [a(t), b(t)],$$
so that, for each fixed $ \zeta$, the map  $t\mapsto y(t, \zeta)$ traces a particle
trajectory.
We denote by $\eta=\eta(t, \zeta)$
the linear density and by
\bel{91}\bfv(t, \zeta)~=~
y_t(t, \zeta)\eeq
the velocity of particles along the curve.
The total mass contained  on a portion of the
curve is thus
\bel{mm}m\Big(\{y(t, \zeta)\,;~ \zeta\in [ \zeta_1,  \zeta_2]\}\Big)~=~\int_{ \zeta_1}^{ \zeta_2} \eta(t, \zeta)\, d \zeta.\eeq

To derive an evolution equation for  $\eta,\bfv$, let
$$\bfv^+(t,x),\quad \bfv^-(t,x),\qquad \rho^+(t,x),\quad \rho^-(t,x),$$
be respectively the velocity and the density (w.r.t.~2-dimensional Lebesgue measure)
of the gas, to the right and to the left of the curve.

\begin{figure}[ht]
\centerline{\hbox{\includegraphics[width=9cm] {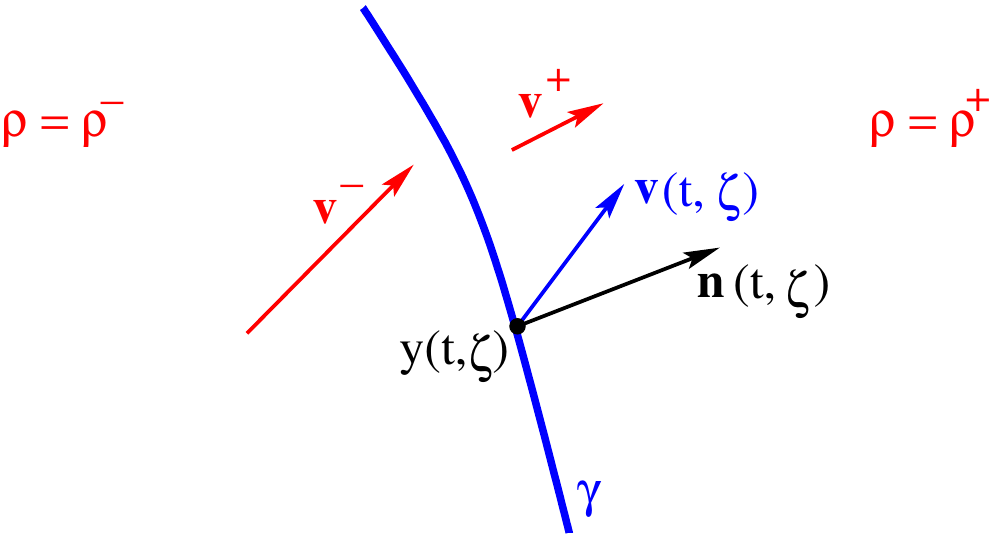}}}
\caption{\small  The variables in the equations (\ref{92}), (\ref{94}), determining the
evolution of a singular curve $\gamma$.}
\label{f:gen77}
\end{figure}

Calling $\bfn(t, \zeta)$ the unit normal vector to the curve at the point $y(t, \zeta)$,
as in Fig.~\ref{f:gen77},
we assume that the inner products satisfy the admissibility condition
\bel{inp}\Big\langle \bfn(t, \zeta),\, \bfv^+(t, y(t, \zeta))\Big\rangle ~\leq~\Big\langle  \bfn(t, \zeta),\, \bfv(t, \zeta)\Big\rangle
~\leq~\Big\langle  \bfn(t, \zeta),\, \bfv^-(t, y(t, \zeta))\Big\rangle .\eeq
These inequalities  imply that particles impinge on the singular curve from both sides
(and then stick to it).

It will be convenient to write the next equation using the
cross product of two vectors $\bfu= (u_1, u_2)$, $\bfv = (v_1, v_2)$ in $\R^2$:
$$\bfu\times \bfv~=~u_1 v_2-u_2 v_1~=~\la \bfu^\perp,\bfv\ra,\qquad \bfu^\perp = (-u_2, u_1).$$

Under the assumption (\ref{inp}), for any $ \zeta_1< \zeta_2$  the rate at which
the mass along the curve  increases is computed by
$$ {d\over dt}\int_{ \zeta_1}^{ \zeta_2} \eta(t, \zeta)\, d \zeta~=~\int_{ \zeta_1}^{ \zeta_2}
\rho^+(\bfv-\bfv^+)\times
y_{ \zeta}\, d \zeta+
\int_{ \zeta_1}^{ \zeta_2}\rho^-(\bfv^--\bfv)\times y_\zeta\, d \zeta \,.$$
Since the above identity holds for all $ \zeta_1< \zeta_2$, this implies
\bel{92}\eta_t(t, \zeta)~=~\rho^+(\bfv-\bfv^+)\times
y_{ \zeta}+\rho^-(\bfv^--\bfv)\times y_\zeta\, ,
\eeq
where the functions $\rho^\pm, \bfv^\pm$
 are evaluated at $(t,y(t, \zeta))$.

A similar argument shows that, by conservation of momentum, one has
\bel{93}
\bigl[\eta(t, \zeta)\, \bfv(t, \zeta)\bigr]_t ~=~[  (\bfv-\bfv^+)\times y_\zeta] \,\rho^+\bfv^++
[  (\bfv^--\bfv)\times y_\zeta] \,\rho^-\bfv^-.\eeq
Combining (\ref{92})  with (\ref{93}), one obtains an
expression for  $\bfv_t(t, \zeta)$, describing the acceleration of a particle along the singular curve.
\bel{94}\bega{rl} \ds y_{tt} (t,\zeta)&=~\bfv_t(t, \zeta)~\ds=~\left( {\eta\bfv\over\eta}\right)_t~=~{(\eta\bfv)_t\over\eta} - {\eta_t\over \eta}\, \bfv
\\[4mm] &\ds=~{1\over\eta} \bigg(
\bigl[  (\bfv-\bfv^+)\times y_\zeta\bigr]\, \rho^+(\bfv^+-\bfv)+
 \bigl[  (\bfv^--\bfv)\times y_\zeta\bigr]\, \rho^-(\bfv^--\bfv)\bigg).
\enda \eeq

%

For  analytic initial data
\bel{icon}y(t_0, \zeta)\,=\,y_0( \zeta),\qquad \eta(t_0, \zeta)\,=\,\eta_0( \zeta),\qquad
\bfv(t_0, \zeta)\,=\,\bfv_0( \zeta),\eeq
a local solution  to the Cauchy problem (\ref{92}), (\ref{94}), (\ref{icon})
was constructed in \cite{BCH}
by an application of the classical Cauchy-Kovalevskaya theorem \cite{Evans, RR, W}.

In the setting considered earlier in \cite{BCH}, the singular curve has already been formed.  Namely, $\eta_0( \zeta)>0$, while $\rho^+, \rho^-$ are bounded analytic functions.
The present situation requires a more careful analysis. Indeed, $\eta_0( \zeta_0)=0$,
while $\rho^-,\rho^+\to +\infty$  as $(t,x)\to (t_0, x_0)$.  In (\ref{94})
these unbounded terms must be partially cancelled by the vanishingly small
terms $\bfv^\pm -\bfv$, so that the right hand side remains integrable.
In Section~\ref{s:5} we will show that this problem can be reformulated as a system of
Briot-Bouquet type, to which the result in \cite{li} can be applied.

 \section{Singularity formation in one space dimension}
 \label{s:3}
\setcounter{equation}{0}

In the one dimensional case, a singularity consists of a point mass.
The formation of a new singularity is similar to the formation of a new shock
in the solution to a scalar conservation law.  Geometric ideas from \cite{Gk2}
can thus be used also in the present situation.

Let $\bar \rho(x)$ and $w(x)$  be the initial density and the initial velocity
of the pressureless gas,
at the point $x\in \R$.
Assuming $w'(x)<0$, the blowup time along the characteristic starting at $x$ is
\bel{blut}\tau(x)~=~{-1\over w'(x)}\,.\eeq
We assume that, at a point $x_0$, the velocity satisfies
\bel{wmin}   w'(x_0) \,=\,-{1\over t_0}\,<\,0,\qquad\quad  w''(x_0)=0,\quad \qquad w'''(x_0)>0\,.\eeq
As a consequence,   $w'$ has a negative minimum at $x_0$, and hence
$t_0 = \tau(x_0)$ is a local minimum of the blow-up
time function $\tau(\cdot)$. In other words, $t_0$
is the first time when the singularity appears.

{}Without loss of generality, to simplify notation we shall assume that $x_0=0$.
To fix ideas,
let the density and velocity functions have the Taylor approximation
\bel{Tar}\bar \rho(x)~=~\sum_{k\geq 0} \rho_k {x^k\over k!}\,
\eeq
\bel{Taw}w(x)~=~\omega_0- \omega_1 x  + \sum_{k\geq 3}  {\omega_k} {x^k\over k!}\,,\eeq
where the coefficients
 satisfy
\bel{cro}
\rho_0>0,\qquad \omega_1>0,\qquad \omega_3 >0.\eeq
In the above setting, we have
\begin{theorem}\label{t:61} Let the initial density and velocity of the 1-dimensional pressureless gas be analytic functions satisfying (\ref{Tar})--(\ref{cro}).
Then there exists a local solution containing a point mass along an analytic curve
$t\mapsto y(t)$, $t\in \,]t_0, t_0+\delta]$. Here $t_0=1/\omega_1$, $y(t_0) = t_0 \omega_0$.

Moreover, the mass has asymptotic  size $m(t)=\sqrt{t-t_0}\cdot \Upsilon(t)$,
where $\Upsilon$ is an analytic function with $\Upsilon(t_0)>0$.
\end{theorem}

\begin{figure}[ht]
\centerline{\hbox{\includegraphics[width=16cm] {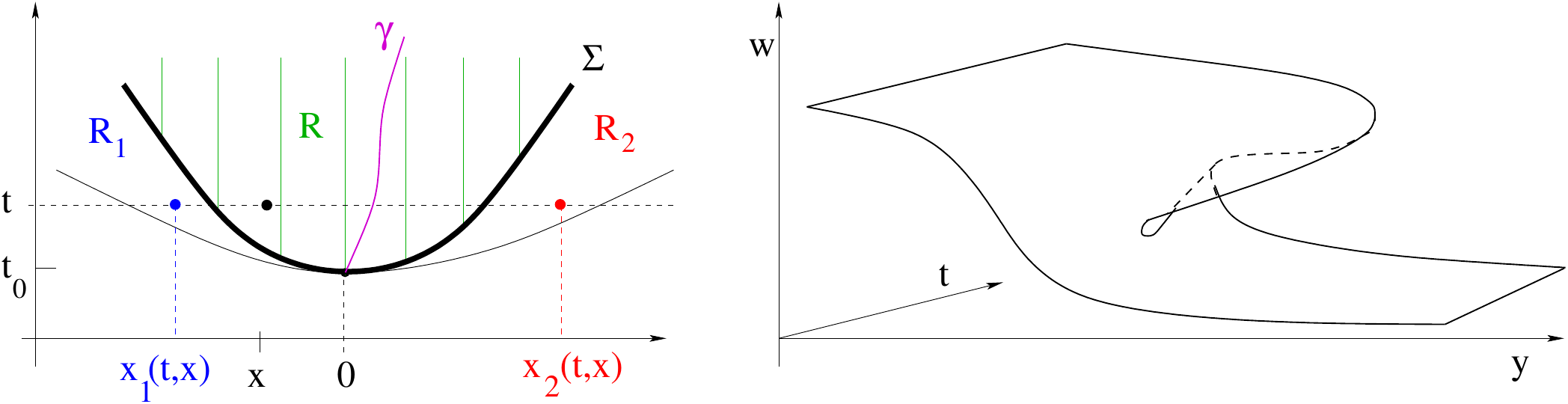}}}
\caption{\small  Left: the regions $R, R_1, R_2$ introduced at (\ref{R0})--(\ref{x34}).
Here $\gamma$ is the stable manifold of the vector field $\bff$ in (\ref{vf}) at the
equilibrium point $(t_0,0)$.
Right: the image of the map $(t,x)\mapsto  (t,y,w)=
\bigl(t, x+ t w(x), w(x)\bigr)$.
For $t>t_0$, this sheet has a fold.   Indeed, if $(t,x)\in R$, there are two additional values
$x_1(t,x)< x < x_2(t,x)$ such that (\ref{x12}) holds. Here $(t, x_1)\in R_1$ while $(t, x_2)\in R_2$.  }
\label{f:gen52}
\end{figure}

\begin{figure}[ht]
\centerline{\hbox{\includegraphics[width=16cm] {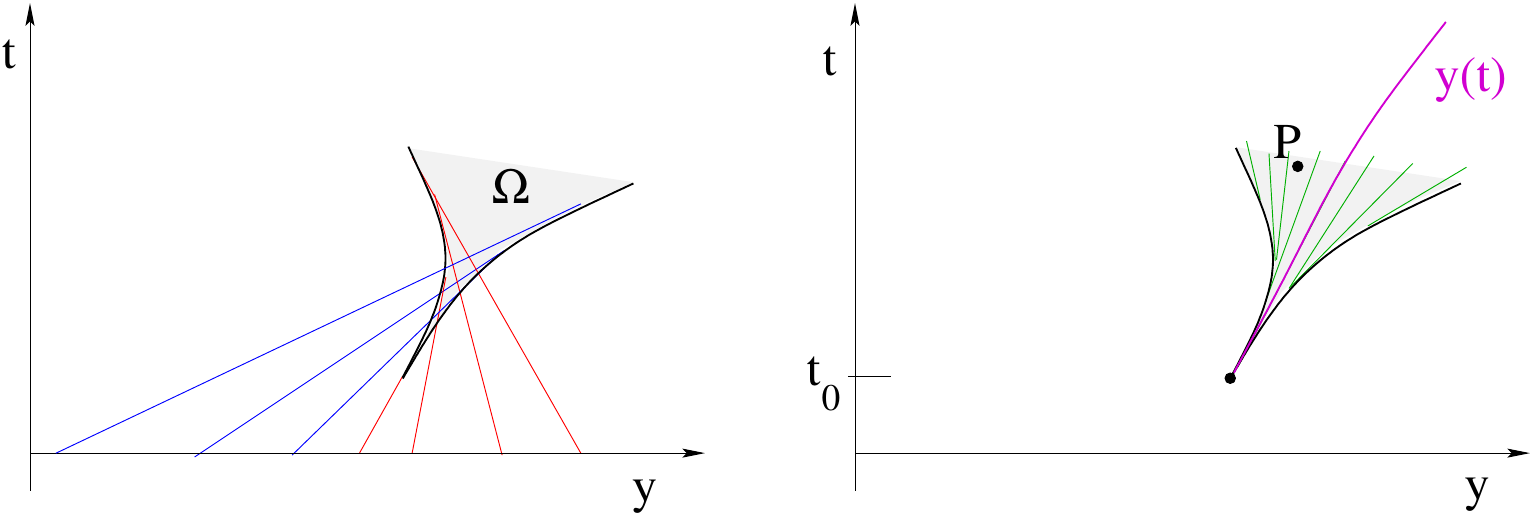}}}
\caption{\small  Every point $P=(t,y)$ in  the region $\Omega$ defined at  (\ref{cuspo})
 is reached by three distinct characteristics. Namely, there exist $x_1<x<x_2$ such that  (\ref{x12}) holds. 
 }
\label{f:gen51}
\end{figure}

{\bf Proof.} {\bf 1.}
To construct the singular curve $y(\cdot)$, following \cite{Gk2} we shall use an alternative
set of coordinates.
 As shown in Fig.~\ref{f:gen52}, left, in the $t$-$x$ plane
we consider the domain
\bel{R0}R~=~\bigl\{ (t,x)\,;~t> \tau(x)\bigr\}~=~\bigl\{(t,x);~1+ t w'(x) <0\bigr\},\eeq
with boundary
\bel{bdR}\Sigma~=~\bigl\{(\tau(x), x)\,;~x\in \R\bigr\}\,.\eeq

 For $(t,x)\in R$, there are two additional points
$(t, x_1(t,x))\in R_1$ and $(t, x_2(t,x))\in R_2$,
with
\bel{x34}x_1(t,x)~<~x~<~x_2(t,x)\,,\qquad\qquad \tau(x_1)\,>\,t,\quad \tau(x)\,<\,t,\quad \tau(x_2)\,>\,t,\eeq
and
such that
\bel{x12}x_1+ t w(x_1)~=~x+ t w(x)~=~x_2+ t w(x_2).\eeq
In other words, the three characteristics starting at $x$, $x_1$ and $x_2$
reach exactly the same point $y$
at time $t$.
As shown in  Figures~\ref{f:gen52} and \ref{f:gen51}, there is a cusp
region covered  by three sheets of the map
$(t,x)\mapsto \bigr(t, x+tw(x), w(x)\bigl)$, namely
\bel{cuspo}\Omega~=~\Big\{ (t,y)\,;~y= x+ t w(x), ~t> \tau(x)\quad\hbox{for some }~x\in\R\Big\}.\eeq
According to (\ref{cuspo}), we shall parameterize points $(t,y)\in \Omega$ by the pull-back
along the characteristic on the middle sheet.
\v
{\bf 2.}
As before, we call $y(t)$ and $m(t)$ the  position and the amount of mass concentrated at the singular point at time $t>t_0$. Writing
\bel{yxt}y(t)~=~x(t) + t w\bigl(x(t)\bigr),\eeq
we seek an evolution equation for $x(\cdot)$ and $m(\cdot)$.
%
%
Conservation of mass and momentum yield
\bel{mint}
m(t)~=~\int_{x_1(t, x(t))}^{x_2(t, x(t))} \bar \rho(  \xi)\, d  \xi,\qquad \qquad
m(t)\dot y(t)~=~\int_{x_1(t, x(t))}^{x_2(t, x(t))} \bar \rho(  \xi)w(  \xi)\, d  \xi\,.\eeq
%
%
%
\v
{\bf 3.}
To proceed, we need an efficient representation of the maps
$x_1(t,x)$, $x_2(t,x)$.
These provide solutions $X=X(t,x)$ to the scalar equation
\bel{XXE}
X+ t w(X) ~=~x+ t w(x).\eeq
Consider  the analytic function
\bel{fde}
f(X, t,x)~\doteq~\left\{ \bega{cl} \ds 1 + t\, {w(X) -w(x)\over X-x}\qquad
&\hbox{if} \quad X\not= x,\\[4mm]
1+ t\, w'(x)&\hbox{if} \quad X= x.\enda
\right.\eeq
Notice that, by replacing the equation (\ref{XXE}) with
\bel{f=0}
f(X,t,x)~=~0,\eeq
we can remove the trivial solution $X=x$.
In view of (\ref{Taw}),  at $t=t_0$  we find
$$f(X, t_0,0) ~=~g(X) \, X^2,$$
where $g$ is an analytic function such that
$g(0)={t_0 \, \omega_3/6}>0$.
According to the Weierstrass Preparation Theorem (see for example \cite{CH} or \cite{GG}), for $(t,x)\approx(t_0,0)$
there exists an analytic function $q$, together with analytic functions $a_1,a_2$ such that
\bel{prep}
 q(X, t,x)\,f(X, t,x) ~=~ X^2 -2 a_1( t,x)  X+ a_0( t, x),\eeq
 with $$a_0(t_0,0)~=~a_1(t_0,0) ~=~ 0,\qquad\qquad q(0,t_0,0)\not= 0.$$
 Since our goal is to find solutions of (\ref{f=0}), the above representation
clarifies the structure of these solutions. Indeed,
 \bel{x12sol}\bega{l}\ds
 x_1(t,x)  ~=~a_1(t,x) - \sqrt{a_1^2(t,x) -  a_0(t,x)}\,,\\[2mm]
 \ds
x_2(t,x) ~=~a_1(t,x) + \sqrt{a_1^2(t,x) - a_0(t,x)}\,.
 \enda
 \eeq

\v
{\bf 4.}
Next,  for any smooth function $\phi$ we can write
\bel{odd}\int_{x-s}^{x+s} \phi(y)\, dy~=~\Phi(x,s^2)\,s\,,\eeq
for a suitable function $\Phi$.  Indeed, the left hand side is an odd function of $s$.
Hence it can be written as a product of $s$ and an even function of $s$.
In connection with (\ref{mint}), (\ref{x12sol}), this yields
\bel{mi2}
m(t)~=~\int_{x_1(t, x(t))}^{x_2(t, x(t))} \bar \rho(\xi)\, d\xi~=~M(t,x) \cdot \sqrt{ a_1^2(t,x) - a_0(t,x)},\eeq
\bel{mi3}
m(t)\dot y(t)~=~\int_{x_1(t, x(t))}^{x_2(t, x(t))} \bar \rho(\xi)w(\xi)\, d\xi~=~
G(t,x) \cdot \sqrt{ a_1^2(t,x) - a_0(t,x)}
\,,\eeq
for suitable analytic functions $M,G$.

Since
$$M(t_0, 0)~=~2\bar\rho(0)~\not= 0,
\qquad\qquad G(t_0, 0) ~=~2\bar\rho(0)\, w(0),$$
 the function
\bel{Fde} F(t,x)~\doteq~{G(t,x)\over M(t,x)} - w(x)\eeq
is  analytic in a neighborhood of $(t_0, 0)$.
By (\ref{mi2})-(\ref{mi3}) we have to solve
$$\dot y~=~\dot x + w(x) + t w'(x) \, \dot x~=~{G(t,x)\over M(t,x)} \,.$$
In view of (\ref{Fde}),
this leads to
\bel{eqx}
(1+t w'(x)) \dot x~=~F(t,x).\eeq

%

Observing that
$$1+ t_0 w'(0)~=~0~=~F(t_0, 0),\qquad\qquad t_0 \,=\, {-1\over w'(0)}\,=\,
{1\over\omega_1}\,,$$
a solution to (\ref{eqx}) can be constructed
by finding an integral curve $\gamma$ of the vector field with analytic  components
\bel{vf}\bff~=~\Big( 1+ tw'(x)\,,~F(t,x)\Big)\eeq
that goes through the equilibrium point $(t_0, 0)$.
We claim that
 this equilibrium is a saddle.

 Indeed, at the point $(t,x)=(t_0, 0)$ the Jacobian matrix has the form
\bel{Jacf}D\bff(t_0,0)~=~
\begin{pmatrix} -\omega_1 & &&0\\[2mm]
\ds{3\omega_1^3\omega_4\over 4\omega_3^2} -
{2\rho_1\omega_1^3 \over\rho_0 \omega_3}
  && &\ds  {3\over 2}  \omega_1 \end{pmatrix}.\eeq
This claim will be proved in Lemma~\ref{l:31} below.
\v
{\bf 5.} Assuming that (\ref{Jacf}) holds, we
can proceed with the proof of the theorem.
Consider the stable manifold for the vector field (\ref{vf}) through
the equilibrium point $(t_0,0)$.
Calling
$t\mapsto x(t)$  the portion of this
manifold for $t>t_0$, we obtain the desired integral curve of $\bff$ through $(t_0,0)$.
In turn, the formula (\ref{yxt})
yields the position of the singularity $y(t)$. We observe that, since the functions $\bar \rho,w$ are analytic, the stable
manifold is analytic as well. See \cite{U} for a general result in this direction.

Finally,
the amount of mass concentrated at $y(t)$ is computed by (\ref{mint}).
More precisely, recalling (\ref{Tar}), by (\ref{mint}) and (\ref{x12sol}) one obtains
\bel{mtser}\bega{rl}
m(t)&\ds=~\sum_{k=0}^\infty \int_{a_1(t, x(t))-\sqrt{ a_1^2(t, x(t)) - a_0(t, x(t))}}^{a_1(t, x(t))+
\sqrt{ a_1^2(t, x(t)) - a_0(t, x(t))}} {\rho_k x^{k}\over k!}\, dx\\[4mm]
&=~\ds \sum_{k=0}^\infty {\rho_k\over (k+1)!} \left[ \sum_{0\leq j\leq k+1\atop j ~{\rm odd}}
2 \Big(\begin{matrix} j\cr k+1\end{matrix} \Big)\big( a_1^2(t, x(t)) - a_0(t, x(t))\bigr)^{j/2}
\cdot a_1(t, x(t))^{k+1-j}\right]\\[4mm]
&=~\sqrt{ a_1^2(t, x(t)) - a_0(t, x(t))}\cdot \Tilde \Upsilon(t)\\[4mm]
&=~\sqrt{t-t_0}\cdot\Upsilon(t),
\enda
\eeq
where $\Tilde \Upsilon, \Upsilon$ are analytic functions.
This completes the proof.
 \endproof

\begin{lemma}\label{l:31} Assume that the functions $\bar\rho, w$ have the Taylor approximations
(\ref{Tar})-(\ref{Taw}).
Then, at the equilibrium point $(t_0,0)$ the Jacobian matrix of $\bff$ is
given by (\ref{Jacf}).
\end{lemma}
{\bf Proof.} 
{\bf 1.}
To compute the partial derivatives of $F$, we begin by observing that, by
(\ref{Taw}),
the equation (\ref{XXE}) can be written as
\bel{xe0} (X-x)+t\left(- \omega_1 (X-x)  +  \sum_{k\geq 3} \omega_k {X^k-x^k\over k!}\right)~=~0.\eeq
Dividing by $X-x$ we remove the trivial solution $X=x$. The remaining two solutions are found by solving
\bel{Xeq2}
1- t\omega_1 + t\sum_{k\geq 3} {\omega_k\over k!}\left( \sum_{j=0}^{k-1} X^j x^{k-1-j}\right)~=~0.\eeq
Collecting lower order terms and
recalling that $t_0\omega_1 = 1$, we obtain
\bel{xe1} -\omega_1 (t-t_0)+
{t\omega_3\over 6} \bigl(X^2 + Xx+x^2\bigr) +  {t\omega_4\over 24
} \bigl(X^3 + X^2x+Xx^2+x^3\bigr)+t \,\Psi(X,x) ~=~0,
\eeq
where $\Psi$ is an analytic function such that
\bel{Psis}\Psi(X,x)~=~\O(1)\cdot \bigl( |X|^4 + |x|^4\bigr).\eeq
Dividing all terms by $t\omega_3/6$ and recalling that $t_0\omega_1 =1$, we are led to the equation
\bel{xe2}f(X,t,x)~\doteq~ X^2 + x X +x^2- {6\omega_1(t-t_0)\over t\omega_3}+{\omega_4\over 4\omega_3}(X^3 +X^2 x + Xx^2 + x^3)+ {6\Psi(X,x)\over \omega_3}~=~0\,.\eeq
Assuming $x=\O(1)\cdot (t-t_0)$, to leading order the solution is given by
\bel{Xlos}
X^\pm(t,x)~=~ \pm \sqrt{{6\omega_1 (t-t_0)\over t \omega_3}}
+\O(1)\cdot (t-t_0)\,.\eeq
Observing that
$$(X^\pm)^3~=~\pm\left( {6\omega_1 (t-t_0)\over t \omega_3}\right)^{3/2}+\O(1)\cdot (t-t_0)^2$$
and inserting this expression in (\ref{xe2}),
to the next order one finds
\bel{X2o}\bega{rl}
X^\pm(t,x)&\approx\,\ds -{x\over 2} \pm \sqrt{{6\omega_1 (t-t_0)\over t \omega_3} \mp{\omega_4\over 4\omega_3} \left( {6\omega_1 (t-t_0)\over t \omega_3}\right)^{3/2} }+\O(1)\cdot (t-t_0)^{3/2}\\[4mm]
&\ds =\,-{x\over 2} \pm \sqrt{{6\omega_1 (t-t_0)\over t \omega_3}} \cdot\sqrt{ 1\mp{\omega_4\over 4\omega_3} \left( {6\omega_1 (t-t_0)\over t \omega_3}\right)^{1/2} } + \O(1)\cdot (t-t_0)^{3/2}
\\[4mm]
&\ds =\,-{x\over 2} \pm \sqrt{{6\omega_1 (t-t_0)\over t \omega_3}}
- {3\omega_1\omega_4 (t-t_0)\over 4 t \omega^2_3} + \O(1)\cdot (t-t_0)^{3/2}
\\[4mm]
&\ds =\,-{x\over 2} \pm \sqrt{{6\omega_1^2 (t-t_0)\over \omega_3}}
- {3\omega_1^2\omega_4 (t-t_0)\over 4  \omega^2_3} + \O(1)\cdot (t-t_0)^{3/2}\\[4mm]
&=\ds~\pm B-C + \O(1)\cdot (t-t_0)^{3/2},
\enda
\eeq
where we set
\bel{BCdef}B~\doteq~\sqrt{{6\omega_1^2 (t-t_0)\over \omega_3}}, \qquad\qquad C~\doteq~{x\over 2} + {3\omega_1^2 \omega_4 (t-t_0)\over 4\omega_3^2}\,.\eeq
\v
{\bf 2.}
Computing a Taylor expansion of the integrand in (\ref{mi2}), for some analytic function
$R_4(\cdot)$ by (\ref{Tar})  we obtain
\bel{mi5}\bega{rl} m(t)&\ds =~\bigg[\rho_0 x+ \rho_1 {x^2\over 2} + \rho_2{x^3\over 3!}+ R_4(x) x^4 \bigg]_{X^-}^{X^+}\\[4mm]
&\ds=~2\rho_0 B - 2\rho_1 BC + {\rho_2\over 3}B^3 + \O(1) \cdot (t-t_0)^2,
\enda
\eeq
A similar expansion holds for (\ref{mi3}). Recalling (\ref{Taw}), for some analytic function $W_4(\cdot)$   we find
\bel{mi6}\bega{rl} m(t)\dot y(t)&\ds=~\bigg[\rho_0\omega_0 x+ (\rho_1\omega_0 - \rho_0\omega_1)  {x^2\over 2} +(\rho_2 \omega_0 - 2\rho_1\omega_1) {x^3\over 6}+
W_4(x) x^4 \bigg]_{X^-}^{X^+}
\\[4mm]
&\ds =~2\rho_0\omega_0 B -2(\rho_1\omega_0 - \rho_0 \omega_1)  B C + {\rho_2\omega_0 - 2\rho_1\omega_1\over 3} B^3 + \O(1)\cdot (t-t_0)^2.
\enda
\eeq
Combining (\ref{mi5}), (\ref{mi6}) and recalling (\ref{Fde}) we finally obtain
{\small
\bel{F1} \bega{rl} F(t,x)&\ds =~{2\rho_0\omega_0 B -2(\rho_1\omega_0 - \rho_0 \omega_1)  B C + (\rho_2\omega_0 - 2\rho_1\omega_1) (B^3/3) \over
2\rho_0 B -2 \rho_1 BC + \rho_2 (B^3/3)} - \omega_0 + \omega_1 x +\O(1)\cdot (t-t_0)^{3/2}
\\[4mm]
&\ds  =~
{\ds \omega_0  -{\rho_1\omega_0 - \rho_0 \omega_1\over \rho_0} C + {\rho_2\omega_0 - 2\rho_1\omega_1\over 6\rho_0} B^2 \over \ds
1- {\rho_1\over \rho_0} C + {\rho_2\over 6 \rho_0}B^2 } - \omega_0 + \omega_1 x +\O(1)\cdot (t-t_0)^{3/2}
\\[4mm]
&=~\ds  \left( \omega_0  -{\rho_1\omega_0 - \rho_0 \omega_1\over \rho_0}
 C + {\rho_2\omega_0 - 2\rho_1\omega_1\over 6\rho_0} B^2 \right)
\left(1+ {\rho_1\over \rho_0} C - {\rho_2\over 6 \rho_0}B^2
\right)\\[4mm]
&\qquad  - \omega_0 + \omega_1 x +\O(1)\cdot (t-t_0)^{3/2}\\[4mm]
&=~\ds
 \omega_1 C
-{\rho_1\omega_1\over 3\rho_0} B^2+ \omega_1 x +\O(1)\cdot (t-t_0)^{3/2}.
\enda
\eeq
}
By (\ref{BCdef}), this yields
\bel{Flinap}
\bega{rl}
F(t,x)&\ds=~ \omega_1 \left({x\over 2} + {3\omega_1^2 \omega_4 (t-t_0)\over 4\omega_3^2}\right)
-{\rho_1\omega_1\over 3\rho_0} \left({6\omega_1^2 (t-t_0)\over \omega_3}
\right)+ \omega_1 x +\O(1)\cdot (t-t_0)^{3/2}\\[4mm]
&\ds=~ \left( {3\omega_1^3 \omega_4\over 4\omega_3^2}
-{2\rho_1\omega_1^2\over \rho_0\omega_3 }\right)(t-t_0) + {3\over 2}\omega_1 x+\O(1)\cdot (t-t_0)^{3/2}.
\enda
\eeq

In view of (\ref{Flinap}), since $w''(0)=0$,
the partial derivatives of the vector field $\bff$ in (\ref{vf}) at the point $(t_0,0)$
are now computed according to (\ref{Jacf}).
\endproof
\v
\begin{remark} {\rm In view of (\ref{Jacf}), at the equilibrium point $(t_0,0)$
the first eigenvalue and eigenvector of the matrix $D\bff$ are computed by
\bel{ee1}\lambda_1\,=\,- \omega_1,\qquad\qquad \bfw_1~=~\left(\bega{c} 1\cr S_0\enda \right),\qquad S_0~=~ {2\over 5\omega_1} \left(
{2\rho_1\omega_1^3 \over\rho_0 \omega_3}  -{3\omega_1^3\omega_4\over 4\omega_3^2} \right).\eeq
The second pair is
\bel{ee2} \lambda_2~=~{3\over 2} \omega_1\, ,\qquad\qquad \bfw_2~=~\left(\bega{c} 0\cr 1\enda\right).\eeq
Recalling that $\omega_1>0$, we conclude that the stable manifold $t\mapsto x(t)$
has asymptotic expansion
\bel{c0def}
x(t)~=~S_0(t-t_0) + \O(1)\cdot (t-t_0)^2,\qquad\qquad
S_0\,=\,{4\over 5}\,{\rho_1\omega_1^2 \over  \rho_0 \omega_3}  -{3\over 10}\,{\omega_1^2\omega_4\over \omega_3^2}.\eeq
}
\end{remark}

\subsection{The second order system.}
In one space dimension, the conservation of  mass and momentum yield the identities  (\ref{mint}).
However, similar identities are not available in dimension $n\geq 2$.
As a preliminary to the analysis
of the 2-dimensional case, we study here the corresponding second order system
of equations.
Denoting derivatives w.r.t.~time with an upper dot,
these take the form
\bel{reds}
\left\{\bega{rl}\dot m&=~(v^--v)\rho^-  -(v^+-v)\rho^+,\\[1mm]
\dot y&=~v,\\[1mm]
\dot v&=~\ds  {\rho^-\over m} (v^--v)^2 -  {\rho^+\over m} (v^+-v)^2 .\enda\right.\eeq
The initial  data at $t=t_0$ are
\bel{bda}
\left\{\bega{rl}m(t_0) & =~0,\\[1mm]
 y(t_0) & =~ x_0+ t_0\,w(x_0),\\[1mm]
v(t_0)&=~w(x_0).\enda\right.\eeq
%
%
Differentiating w.r.t.~time the identity
\bel{y0} y~=~x+t w(x),\eeq
 one finds
\bel{y1}\dot y~=~w(x) + \bigl[ 1+ t w'(x)\bigr] \dot x,\eeq
\bel{y2}\ddot y~=~\bigl[ 1+ t w'(x)\bigr]  \ddot x  +2 w'(x) \dot x + t\,w'' (x) \dot x^2. \eeq
In the following, we denote the density and speed behind and ahead of the singularity as functions
of the variables $t,x$, namely
\bel{rvpm}
\rho^\pm(t,x)~=~ \bigl[ 1+ t w'\bigl(X^\pm(t,x)\bigr)\bigr]^{-1} \bar \rho\bigl(X^\pm(t,x)\bigr),
\qquad\qquad v^-(t,x)~=~w\bigl(X^-(t,x)\bigr).\eeq

where $X^+, X^-$ are the two solutions of (\ref{xe0}), computed at (\ref{X2o})-({\ref{BCdef}).

Setting $S= \dot x$, in these new variables
the system (\ref{reds}) takes the form
\bel{SO1}
\left\{ \bega{rl} \dot m&=~\ds\Big(w(x) + \bigl[ 1+ t w'(x)\bigr] S- v^+\Big)\rho^+
- \Big(w(x) + \bigl[ 1+ t w'(x)\bigr] S- v^-\Big)\,\rho^-,\\[4mm]
\dot x&=~S,\\[3mm]
\dot S&=~\ds   \bigl[ 1+ t w'(x)\bigr]^{-1} \Bigg\{  \Big(w(x) +\bigl[ 1+ t w'(x)\bigr] S- v^-\Big)^2{\rho^-\over m}
\\[4mm]
&\qquad\qquad\ds - \Big(w(x) + \bigl[ 1+ t w'(x)\bigr]  S- v^+\Big)^2 {\rho^+\over m}
-2 w'(x) S - t\,w'' (x) S^2\Bigg\}.
\enda
\right.\eeq
Of course, the solution to (\ref{SO1}) is the same one described in Theorem~\ref{t:61}.
Toward a future extension to the 2-dimensional case we remark that,
while (\ref{SO1}) contains three equations,
only two initial conditions are given, namely
\bel{mx00} m(t_0) = 0,\qquad\qquad x(t_0)=0.\eeq
At this stage, an important  observation is that the initial value  $S_0=\dot x(t_0)$
can be uniquely determined
by imposing the integrability of the right hand side of the last equation in (\ref{SO1}).
Indeed, the map
$t\mapsto \bigl[ 1+ t w'(x)\bigr]^{-1} $ is not integrable
in a neighborhood of $t=t_0$.  To achieve a well defined solution
it is essential that,  on the right hand side of this equation,
the second factor within brackets vanishes at $t=t_0$.

Assume that
\bel{xe} x(t)~=~S_0(t-t_0) + \O(1)\cdot (t-t_0)^2,\eeq
for some constant $S_0$ to be determined.
By (\ref{X2o}) this implies
\bel{XSpm} X^\pm(t)~=~- {S_0(t-t_0)\over 2} \pm \sqrt{{6\omega_1^2 (t-t_0)\over \omega_3}}
- {3\omega_1^2\omega_4 (t-t_0)\over 4  \omega^2_3} + \O(1)\cdot (t-t_0)^{3/2}.\eeq
Recalling that $\omega_1 t_0=1$, from (\ref{Tar})-(\ref{Taw}) we obtain
\bel{oom}
1+tw'\bigl(x(t)\bigr)~=~-\omega_1(t-t_0) + \O(1)\cdot (t-t_0)^2,\eeq
$$v^--v^+~=~w\bigl( X^-(t)\bigr) - w\bigl( X^+(t)\bigr)~=~2\omega_1\sqrt{ 6\omega_1^2 (t-t_0)\over \omega_3} + \O(1)\cdot (t-t_0).$$
$$v^--w(x)~=~\omega_1^2\sqrt{ 6 (t-t_0)\over \omega_3} + \O(1)\cdot (t-t_0).\,,$$
In addition, by (\ref{XSpm}) we have the
leading order approximations
\bel{denn}\bega{rl}
1+ t w'(X^\pm)&\ds=\,-\omega_1 (t-t_0) + t_0 \omega_3 {(X^\pm)^2\over 2}+t_0\omega_4{(X^\pm)^3\over 6}
+\O(1)\cdot (t-t_0)^2\\[4mm]
&=~
2\omega_1(t-t_0)\bigl(1\mp E\cdot \sqrt{t-t_0}\,\bigr)+\O(1)\cdot (t-t_0)^2\,,\enda
\eeq
where
\bel{Edef}
E\,\doteq\,\frac{\omega_3^2S_0-\frac12\omega_1^2\omega_4}{4\omega_1\omega_3}\sqrt{\frac{6}{\omega_3}}.\eeq
In turn, this yields
\bel{rh1}\rho^\pm~=~\frac{\bar\rho(X^\pm)}{1+tw'(X^\pm)}~=~\frac{1}{2\omega_1(t-t_0)}\left[\rho_0\pm\left(\rho_1\omega_1\sqrt{\frac{6}{\omega_3}}
+\rho_0 E\right)\sqrt{t-t_0}\right]+\O(1),
\eeq
$$\dot m~=~(v^--v^+) \rho^+ +
 \Big(w(x) + \bigl[ 1+ t w'(x)\bigr] S- v^-\Big)\,(\rho^+-\rho^-)
 ~=\ds~ { \omega_1 \rho_0\over \sqrt{t-t_0} }\, \sqrt{6\over\omega_3}+\O(1)\,,
 $$
\bel{mt}m(t)~=~2\omega_1 \rho_0  \sqrt{6(t-t_0)\over\omega_3}+\O(1)\cdot (t-t_0).\eeq

%
%

In particular, on the right hand side of (\ref{SO1}), the coefficients of the terms containing $S^2=\dot x^2$ have size
\bel{qterms}\left\{\bega{rl} \ds
\bigl[ 1+ tw'(x)\bigr]^2 {\rho^--\rho^+\over m}&\ds=~\O(1)\cdot (t-t_0)^2\cdot {(t-t_0)^{-1/2}\over (t-t_0)^{1/2}}~=~\O(1)\cdot (t-t_0),\\[4mm]
t w''(x)&=~\O(1)\cdot (t-t_0).\enda\right.\eeq
Therefore, even after multiplication by  $\bigl[ 1+tw'(x)\bigr]^{-1}=\O(1)\cdot (t-t_0)^{-1}$,
these terms are integrable.
Thanks to this observation, the integrability of the right hand side in the last equation of (\ref{SO1})
yields an equation which is linear in $S_0$.
%
%
Recalling that $x(t)\approx (t-t_0) S_0$,
we compute the approximations
\bel{xd0}\bega{rl} \dot S&\ds\approx~{-1\over \omega_1(t-t_0)}
\cdot\bigg\{  \Big(\omega_0 -  \omega_1 x(t) -
\omega_1  (t-t_0)S_0  - w\bigl(X^-(t)\bigr)\Big)^2{\rho^-(t)\over m(t)}
 \\[4mm]
&\qquad\qquad\qquad\ds -\Big(\omega_0 -  \omega_1 x(t) -
\omega_1  (t-t_0)S_0
- w\bigl(X^+(t)\bigr)\Big)^2{\rho^+(t)\over m(t)}+
2\omega_1 S_0 \bigg\}\\[5mm]
&\ds\approx~{-1\over \omega_1(t-t_0)}
\cdot\bigg\{  \omega_1^2 \bigl(X^-(t)\bigr)^2{ \rho^-(t)\over m(t)} - \omega_1^2 \bigl(X^+(t)\bigr)^2{ \rho^+(t)\over m(t)} \\[4mm]
&\qquad\qquad\qquad\ds-\omega_1^2\frac{4S_0(t-t_0)}{m(t)}\Big(\rho^-(t)X^-(t)-\rho^+(t)X^+(t)\Big)+
 2 \omega_1 S_0\bigg\}\\[5mm]
&\ds=~{-1\over t-t_0} \cdot\Bigg\{ {\omega_1\over m(t)} \bigg(  (X^-)^2 \bigl[\rho^-- \rho^+\bigr] +  \bigl[(X^-)^2 - (X^+)^2\bigr] \rho^+\\[4mm]
&\qquad\qquad\qquad \qquad\ds
-4S_0(t-t_0)\big(\rho^-X^--\rho^+X^+\big)\bigg)
+2 S_0\Bigg\}
\,.\enda\eeq
We need to impose that the factor within braces
vanishes as $t\to t_0$.   Recalling (\ref{denn}) and (\ref{rh1}),
and then using (\ref{XSpm}),
we find
\bel{S00}\bega{l}\ds
{\omega_1\over m(t)} \Big(  (X^-)^2 \bigl[\rho^-- \rho^+\bigr] +  \bigl[(X^-)^2 - (X^+)^2\bigr] \rho^+
 {-4S_0(t-t_0)\big(\rho^-X^--\rho^+X^+\big)}\Big)
+2 S_0\\[4mm]
\ds \qquad \approx{\omega_1\over m(t)} \Bigg\{-\frac{6\omega_1^2}{\omega_3}(t-t_0)\cdot\left(\rho_1\sqrt{\frac{6}{\omega_3}}+\frac{\rho_0}{\omega_1}E\right){1\over \sqrt{t-t_0}}\\[4mm]
\ds\qquad\qquad
+\frac{2S_0\omega_3^2+3\omega_1^2\omega_4}{\omega_3^2}\omega_1\sqrt{\frac{6}{\omega_3}}(t-t_0)^{3/2}\cdot\frac{1}{2\omega_1(t-t_0)}\left[\rho_0+\left(\rho_1\omega_1\sqrt{\frac{6}{\omega_3}}
+\rho_0 E\right)\sqrt{t-t_0}\right] \\[4mm]
\ds \qquad\qquad -4S_0(t-t_0)\cdot\left({-\rho_0\over 2\omega_1(t-t_0)} 2\omega_1 \sqrt{6(t-t_0)\over\omega_3}\,\right)\Bigg\}+2S_0\,.
\enda
\eeq
Finally, recalling the expressions for $m(t)$ at (\ref{mt}) and for $E$ at (\ref{Edef}),
imposing that the quantity in (\ref{S00}) approaches zero  as $t\to t_0$, we obtain
{\small
\bel{S0lineq}\bega{l}\ds
0~=~{\omega_1\over  2\rho_0\omega_1\sqrt{6/\omega_3}} \left\{-\frac{6\omega_1^2}{\omega_3}\cdot\left(\rho_1\sqrt{\frac{6}{\omega_3}}+\frac{\rho_0}{\omega_1}E\right)
+\frac{2S_0\omega_3^2+3\omega_1^2\omega_4}{\omega_3^2}\sqrt{\frac{6}{\omega_3}}\cdot\frac{\rho_0}{2} + 4S_0\rho_0\sqrt{\frac{6}{\omega_3}}\;\right\}+2S_0\\[5mm]
\ds \quad=~{1\over  2\rho_0} \left\{-\frac{6\omega_1^2}{\omega_3}\cdot\left(\rho_1+\frac{\rho_0}{\omega_1}\cdot\frac{\omega_3^2S_0-\frac12\omega_1^2\omega_4}{4\omega_1\omega_3}\right)
{ +\frac{2S_0\omega_3^2+3\omega_1^2\omega_4}{\omega_3^2}\cdot\frac{\rho_0}{2}} + 4S_0\rho_0\;\right\}+2S_0\\[5mm]
\ds\quad=~\frac{1}{2\rho_0}\left(\frac{7}{2}S_0\rho_0-6\frac{\omega_1^2\rho_1}{\omega_3}+\frac94\cdot\frac{\rho_0\omega_1^2\omega_4}{\omega_3^2}\right)
+2S_0\,.
\enda
\eeq}
This yields
\bel{c00}S_0~=~{4\over 5} {\rho_1\omega_1^2\over \rho_0\omega_3} -
{3\over 10} {\omega_1^2\omega_4\over\omega_3^2}\,.\eeq
Notice that, by this different approach, we again recover the same value for $S_0$ as
in (\ref{c0def}).

 \section{Singularity formation in 2D}
 \label{s:4}
\setcounter{equation}{0}
Beginning with  this section, we study the formation of a new singularity in  two space dimensions.
Assuming that the initial data $\bar \rho, \bfw$ in (\ref{ic}) are analytic, we wish to construct
a solution  of (\ref{SP}) beyond the time $t_0$ when the singularity first appears.
Referring to Fig.~\ref{f:gen60}, for $t>t_0$ we expect that the solution will be smooth outside
a singular curve $\gamma^t$, where a positive amount of mass is concentrated.
Outside $\gamma^t$, the velocity  and the density of the pressureless gas are
still computed by (\ref{vel}) and (\ref{den}) respectively.   To completely describe this
solution, it thus remains to determine the position of the singular curve $\gamma^t$, and
the velocity and density of particles along this curve.

We observe that, in a neighborhood of a point $Q$ in the relative interior of $\gamma^t$, where the density $\eta$ is strictly positive, the singular curve is described by the system of equations
(\ref{92})--(\ref{94}). Its local construction has already been
worked out in \cite{BCH}.
Aim of this and the next section is to construct the singular curve $\gamma^t$ in a neighborhood of an endpoint
$P$.  The same  construction remains valid near the point $P_0$ where
the singularity first appears.

While the density $\eta$ and the velocity $\bfv$ of particles along the singular curve
are still determined by the balance of mass and momentum (\ref{92})-(\ref{93}),
we can no longer rely on the 1-dimensional formulas (\ref{mint}).
Always referring to  Fig.~\ref{f:gen60}, near an endpoint
$P\in\gamma^t\cap \gamma^\sharp$, one has
$$\rho^\pm\to+\infty,\qquad\qquad \eta\to 0.$$
Therefore, both equations (\ref{92}) and (\ref{94}) contain unbounded terms.
Toward the analysis of these two equations, the key step is to get a
manageable expression for the  functions $\bfv^\pm$, $\rho^\pm$, measuring the
velocity and the density of the gas in a neighborhood of the singular curve.
Of course, $\rho$ and $\bfv$ can be determined a priori, in terms of
(\ref{vel})-(\ref{den}).     However, in the region where the singular curve is located,
 every point $(t,y)$ can be reached by three distinct characteristics:
$$y~=~X_1 + t \bfw(X_1) ~=~X_2 + t \bfw(X_2)~=~X_3 + t \bfw(X_3).$$
The points $X_1,X_2,X_3\in \R^2$ are thus obtained as the roots of a cubic-like equation.
This would lead
to very messy computations.  In alternative,
following the approach used in
Section~\ref{s:3} for the 1-dimensional
problem, we use one of these solutions as independent coordinates $(t, x)$, and denote by
$X_1(t,x)$, $X_2(t,x)$ the two remaining solutions of
$$X+ t \bfw(X)~=~x + t \bfw(x).$$
By doing so, the functions $X_1, X_2$ are obtained by solving an equation
of quadratic type, and can be conveniently expressed in terms of square roots.
As in \cite{Gk2}, here $x$ is the initial point of the characteristic reaching $y= x+t\bfw(x)$ along the ``middle sheet".
It can be distinguished from the other two characteristics by the relations
$$ \tau(X_1) > t,\qquad \tau(x) < t,\qquad \tau(X_2) > t.$$
Here and in the sequel, $\tau(x)$ denotes the first time $t>0$ where the matrix
$I + t D\bfw(x)$ is not invertible.

Ultimately, the goal of this section is to obtain a local asymptotic expansion for the
velocity and the density of the gas:
$$\bfv^-(t,x)~=~\bfw\bigl( X_1(t,x)\bigr), \qquad \bfv^+(t,x)~=~\bfw\bigl( X_2(t,x)\bigr),$$
$$\rho^-(t,x)~=~{\bar \rho\bigl( X_1(t,x)\bigr)\over
\det \bigl( I + t D\bfw (X_1(t,x))\bigr)}\,,\qquad\rho^+(t,x)~=~{\bar \rho\bigl( X_2(t,x)\bigr)\over
\det \bigl( I + t D\bfw (X_2(t,x))\bigr)}\,.$$
This will be achieved in Lemma~\ref{l:42}.

We now start working out details. Let $\bfw$ be
 a vector field satisfying the assumptions in {\bf (A1)-(A2)}.
 In particular,  the
$2\times 2$ Hessian
matrix $D^2\lambda_1(x_0)$ of second order derivatives of $\lambda_1$,
computed at $x=x_0$, is strictly positive definite.  Hence the same is true of
the Hessian matrix $D^2\tau(x_0)$.

We denote by
$\bigl\{r_1(x),\,r_2(x)\bigr\}$ 
a basis of right eigenvectors
of $D\bfw(x)$, normalized so that
\bel{eigv} \bigl|r_1(x)\bigr|~=~\bigl|r_2(x)\bigr|~=~1.\eeq

\begin{figure}[ht]
\centerline{\hbox{\includegraphics[width=14cm] {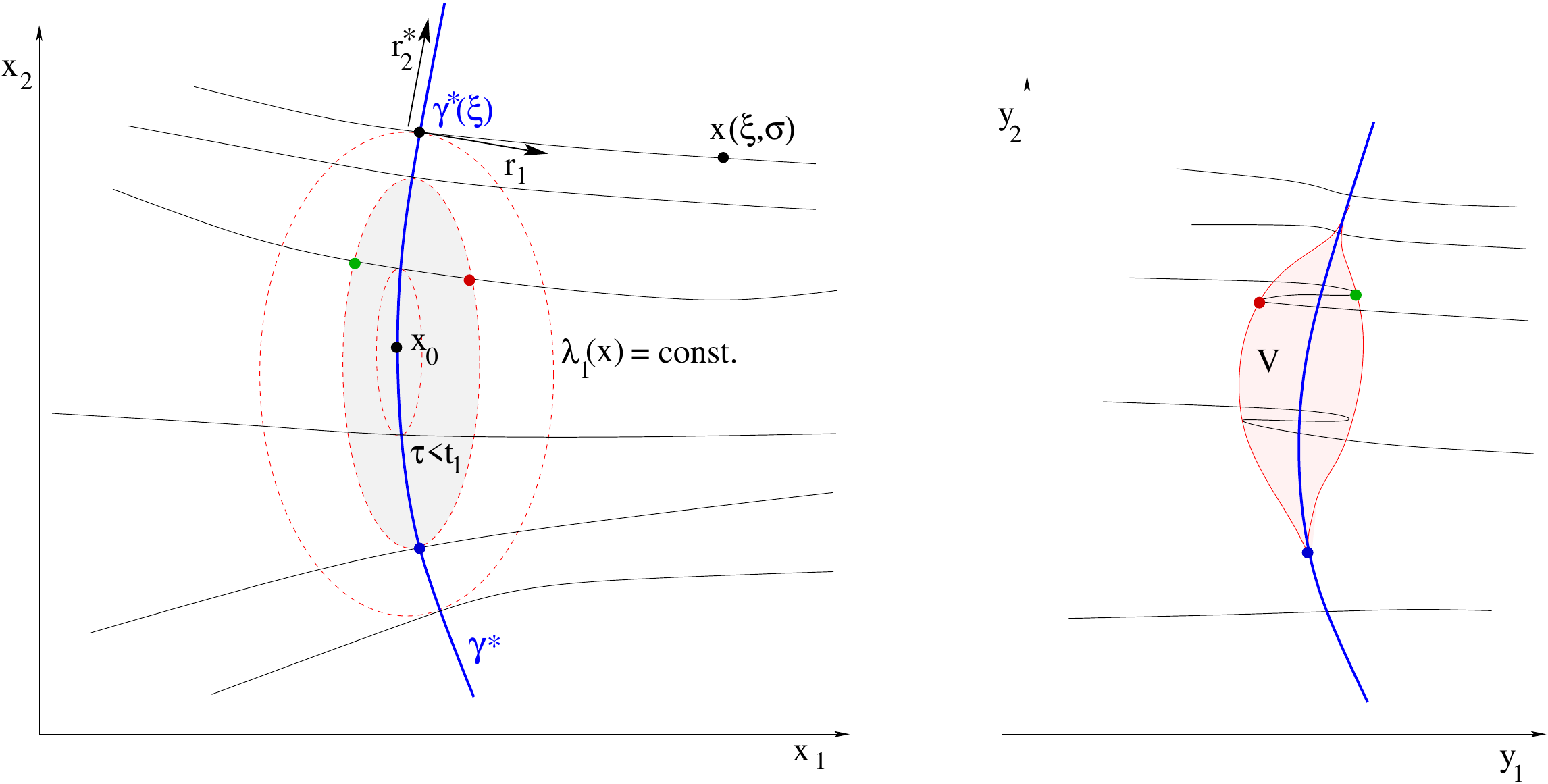}}}
\caption{\small Left: the integral curves of the corresponding eigenvector $r_1$, and the level sets of the eigenvalue $\lambda_1(x)$ of $D\bfw(x)$,
which coincide with the level curves of the blow-up time
$\tau(x) = -1/\lambda_1(x)$.
Here $\gamma^*$
is the set of points where these two families of curves are tangent.
Moreover, $x(\xi,\sigma)$ is the point reached by starting at $\gamma^*(\xi)$ and moving along the integral curve of the vector field $r_1$ for a time $\sigma$.
Right: the shaded region $V$
is the image of the set $\{x;\,\tau(x)<t_1\}$, under the map
$x\mapsto x+ t_1\bfw(x)$.  Each point $y\in V$ has three pre-images. }
\label{f:gen42}
\end{figure}

In the above setting, by the implicit function theorem there exists a smooth curve
\bel{gds}\gamma^*~\doteq~\Big\{ x\in\R^2\,;~\la\nabla
\lambda_1(x)\,,\, r_1(x)\ra \,=\,0\Big\},\eeq
defined near  $x_0$ (see Fig.~\ref{f:gen42}, left).
Indeed, since $\nabla \lambda_1(x_0)=0$ and $D^2\lambda_1(x_0)$ is
positive definite, we have
$$\nabla\la\nabla
\lambda_1(x)\,,\, r_1(x)\ra\bigg|_{x=x_0}~=~\bigl(D^2\lambda_1(x_0)
\bigr)\, r_1(x_0)~\not=~ 0.$$
We parameterize this curve by arc-length:
$\xi\mapsto \gamma^*(\xi)$, so that $\gamma^*(0)= x_0$.
Moreover,
on a neighborhood of
$x_0$, we choose a system of coordinates $\xi,\sigma$
so that (see Fig.~\ref{f:gen42}, left)
\bel{xxs}x(\xi,\sigma)~=~
\exp(\sigma r_1)(\gamma^*(\xi)).
\eeq
Here $\sigma\mapsto \exp(\sigma r_1)(\bar x)$
denotes the solution to the ODE
$$\dot x~=~r_1(x),\qquad\qquad  x(0)=\bar x.$$
To simplify our notation, for $i=1,2$ we shall write
\bel{defs}\left\{\bega{l}\lambda_i(\xi)\,=\,\lambda_i(\gamma^*(\xi)),\\[3mm]
r_i(\xi)\,=\,r_i(\gamma^*(\xi)),\\[3mm]
\tau(\xi)\,=\,\tau(\gamma^*(\xi)),\\[3mm]
\bar\rho(\xi)\,=\,\bar\rho(\gamma^*(\xi)),\\[3mm]
\bfw(\xi)\,\doteq\, \bfw(\gamma^*(\xi)),\enda\right.
\qquad
\qquad \left\{\bega{l} \lambda_i(\xi,\sigma)\,=\,\lambda_i(x(\xi,\sigma)),\\[3mm]
r_i(\xi,\sigma)\,=\,r_i(x(\xi,\sigma)),
\\[3mm]\tau(\xi,\sigma)\,=\,\tau(x(\xi,\sigma)),\\[3mm]
\bar\rho(\xi,\sigma)\,=\,\bar\rho(x(\xi,\sigma)),\\[3mm]
 \bfw(\xi,\sigma)\,=\,\bfw(x(\xi,\sigma)).\enda\right.\eeq
Notice that, in view of (\ref{xxs}), the above notation implies
$\lambda_i(\xi)=\lambda_i(\xi,0)$, ~$\ldots$~, $\bfw(\xi) = \bfw(\xi,0)$.
In addition, for every $\xi$ we define the unit tangent vector to the curve $\gamma^*$ by setting
\bel{r*}r_2^*(\xi)~\doteq~{d\over d\xi} \gamma^*(\xi).\eeq
We observe that, in general, $r_2^*(\xi)\not= r_2(\xi)$.   However, our assumptions imply that
$\bigl\{ r_1(\xi), r_2^*(\xi)\bigr\}$ is still a basis of $\R^2$.

\begin{figure}[ht]
\centerline{\hbox{\includegraphics[width=11cm] {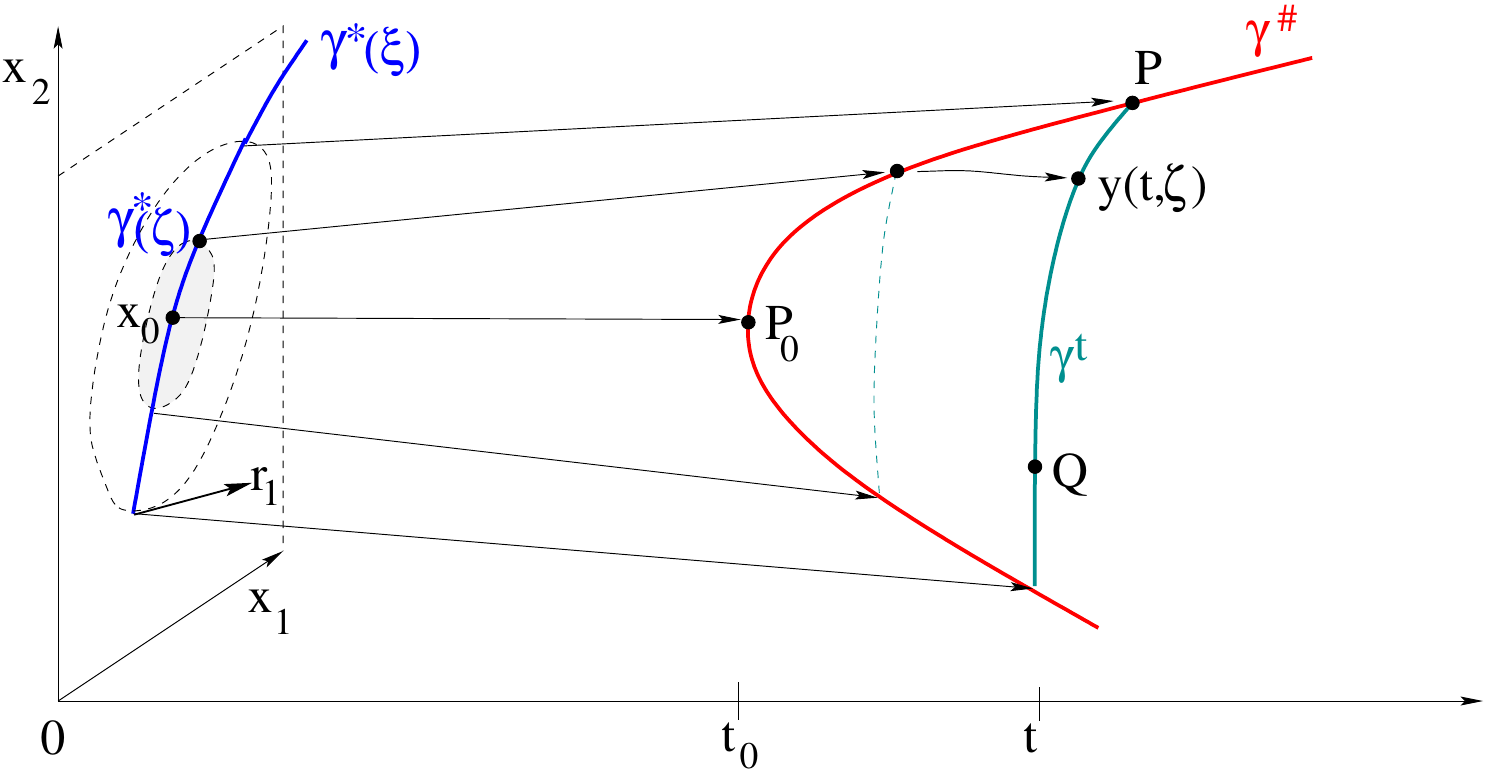}}}
\caption{\small  Here $\gamma^*$ is the curve in the $x_1$-$x_2$ plane, constructed in Fig.~\ref{f:gen42}.  The image of $\gamma^*$
under the map $x \mapsto \bigl(\tau(x), x+ \tau(x) \bfw(x)\bigr)$ is a curve $\gamma^\sharp$
in $t$-$x$ space.
At each time $t>t_0 = \tau (x_0)$,  the singular curve $\gamma^t$ (where a positive amount of mass is concentrated) is parameterized by $\zeta\mapsto y(t,\zeta)$, where $\tau(\zeta)\leq t$.
Here $y(t,\zeta)$ is the position at time $t$ of the particle which was located at
$\gamma^*(\zeta)$
at the initial time $t=0$.
The endpoints of $\gamma^t$ lie on $\gamma^\sharp$.}
\label{f:gen60}
\end{figure}

Similarly to the 1-dimensional case,  for  $t>t_0$ the singular curve $\gamma^t$ will be located
within the region $V$ in the interior of the cusp covered by three sheets of characteristics.
Based on this fact, we define the set
\bel{Rdef2}
R~\doteq~\bigl\{ (t,\xi,\sigma)\,;~\tau(\xi,\sigma)<t\bigr\},\eeq
and its boundary
$$\Sigma~\doteq~\bigl\{ (t,\xi,\sigma)\,;~\tau(\xi,\sigma)= t\bigr\}.$$
For each $(t,\xi,\sigma)\in R$, there exist two other characteristics reaching
the same point at time $t$.   We call
$$\Big(\xi_i(t,\xi,\sigma)\,, ~\sigma_i(t,\xi,\sigma)\Big),\qquad i=1,2,$$
the coordinates of the initial points of these characteristics.
In other words, the points
\bel{xi12}x\bigl(\xi_1(t,\xi,\sigma),\,\sigma_1(t,\xi,\sigma)\bigr),
\qquad
x\bigl(\xi_2(t,\xi,\sigma),\,\sigma_2(t,\xi,\sigma)\bigr),\eeq
are the initial points of the additional two characteristics that at time $t$ reach the same
point
\bel{yxs}y~=~x(\xi,\sigma)+ t\, \bfw\bigl(x(\xi,\sigma)\bigr).\eeq
As a consequence, using the shorter notation introduced at (\ref{defs}), at time $t$
the values of the velocity $\bfv^\pm$  at the point $y$ are given by
\bel{vpm}\bega{rl} \bfv^-(t,y)&=~\bfw\Big(\xi_1(t,\xi,\sigma),\,\sigma_1(t,\xi,\sigma)\Big),\\[3mm]
\bfv^+(t,y)&=~\bfw\Big(\xi_2(t,\xi,\sigma),\,\sigma_2(t,\xi,\sigma)\Big).\enda
\eeq
By (\ref{den}),
similar formulas can be obtained for the densities $\rho^\pm(t,y)$. Namely,
\bel{rpm}\bega{rl} \rho^-(t,y)&\ds
=~\frac{\bar{\rho}\bigl(\xi_1(t,\xi,\sigma),\sigma_1(t,\xi,\sigma)\bigr)}
{\det\Big[I+t\,D_x\bfw\bigl(\xi_1(t,\xi,\sigma),\sigma_1(t,\xi,\sigma)\bigr)\Big]},\\[4mm]
\rho^+(t,y)&\ds=~\frac{\bar{\rho}\bigl(\xi_2(t,\xi,\sigma),\sigma_2(t,\xi,\sigma)\bigr)}
{\det\Big[I+t\,D_x\bfw\bigl(\xi_2(t,\xi,\sigma),\sigma_2(t,\xi,\sigma)\bigr)\Big]}.\enda\eeq
We observe that, for $i=1,2$, the determinant is computed by
\bel{det12}\bega{l}\det\Big[I+t\,D_x\bfw\bigl(\xi_i(t,\xi,\sigma),\sigma_i(t,\xi,\sigma)\bigr)\Big]
\\[3mm]
\qquad =~\Big[1+t\lambda_1\bigl(\xi_i(t,\xi,\sigma),\sigma_i(t,\xi,\sigma)\bigr)\Big]
\Big[1+t\lambda_2\bigl(\xi_i(t,\xi,\sigma),\sigma_i(t,\xi,\sigma)\bigr)\Big].\enda\eeq
The first factor on the right hand side approaches zero as $t\to {-1/\lambda_1}$, while the second one remains uniformly positive.

\v
\subsection{Constructing the points $x_i(t,\xi,\sigma)$.}
Given $(t, \xi,\sigma)$, with $t>\tau(\xi,\sigma)$, we seek a formula describing
the two points in (\ref{xi12}).  For $i=1,2$, the vector equation
$$x(\xi_i, \sigma_i)+ t \, \bfw\bigl( x(\xi_i, \sigma_i)\bigr) ~=~x(\xi,\sigma)+ t\, \bfw\bigl(x(\xi,\sigma)\bigr)$$
is equivalent to two scalar equations, obtained by taking the inner products with $r_1(\xi)$ and  with its
orthogonal vector.   Here and in the sequel, $r_1^\perp(\xi)$ denotes the unit vector obtained by rotating
$r_1(\xi)$ counterclockwise by an angle $\pi/2$.
\bel{le}\left\{
\bega{rl} \Phi_1(t,\xi, \sigma,\xi_i,\sigma_i)&\doteq~\Big\langle r_1(\xi)\,,~ x(\xi_i, \sigma_i)+ t \, \bfw\bigl( x(\xi_i, \sigma_i)\bigr) - x(\xi,\sigma)- t\, \bfw\bigl(x(\xi,\sigma)\bigr)
\Big\rangle~=~0,\\[4mm]
\Phi_2(t,\xi, \sigma,\xi_i,\sigma_i)&\doteq~\Big\langle r_1^\perp(\xi)\,,~ x(\xi_i, \sigma_i)+ t \, \bfw\bigl( x(\xi_i, \sigma_i)\bigr) - x(\xi,\sigma)- t\, \bfw\bigl(x(\xi,\sigma)\bigr)
\Big\rangle~=~0.\enda\right.\eeq
The equations (\ref{le}) have to be solved for the two variables $\xi_i, \sigma_i$.
Trivially,
$$ \Phi_1(t,\xi, \sigma,\xi,\sigma) ~=~\Phi_2(t,\xi, \sigma,\xi,\sigma)~=~0.$$
Of course, we are interested in the two additional nontrivial
solutions.
When $t=\tau(\xi)$, $\xi_i=\xi$  and $\sigma_i=\sigma=0$, we observe that the vector
$r_2^*(\xi)=
{\partial \over\partial \xi}\gamma^*(\xi)$ introduced at (\ref{r*}) is not parallel to $r_1(\xi)$.
Say,
\bel{rc12}r_2^*(\xi)~=~c_1(\xi) r_1(\xi) + c_2(\xi) r_2(\xi),\eeq
with $c_2(\xi)\not= 0$.
Therefore,
\bel{p2xi} \bega{rl} p_0(\xi)&\ds \doteq~\frac{\partial \Phi_2}{\partial \xi_i}\bigl(\tau(\xi), \xi, 0,\xi_i,0\bigr)\bigg|_{\xi_i=\xi}\,=\,\Big\langle r_1^\perp(\xi),
\,\bigl(I+\tau(\xi)D\bfw\bigr)r_2^*(\xi)\Big\rangle\\[3mm]
&=\, \Big\langle r_1^\perp(\xi),
\,\bigl(1+\tau(\xi)\lambda_2(\xi)\bigr)c_2(\xi)r_2(\xi)\Big\rangle\,\neq\, 0.\enda\eeq
The second equation in (\ref{le}) can thus be solved for the variable
$\xi_i$ by the implicit function theorem, say
\bel{xii}\xi_i~=~\phi(t, \xi,\sigma, \sigma_i),\eeq
for some analytic function $\phi$ with
\bel{pze}
\phi\bigl(t, \xi, \sigma,\sigma\bigr)~=~\xi\eeq
for every $t,\xi,\sigma$.  The following computations will yield a leading order expansion for the function $\phi$.

Since the map $\sigma\mapsto \tau(\xi,\sigma)= - 1/\lambda_1(\xi,\sigma)$
has a positive local minimum at $\sigma=0$, a  Taylor expansion yields
\bel{Tlam}\lambda_1(\xi,\sigma)~=~-\omega_1(\xi) + \sum_{k\geq 3} \omega_k(\xi) {\sigma^{k-1}\over (k-1)!}
\eeq
for some analytic functions $\omega_j(\xi)$, with
$$\omega_1(\xi)\,>\,0,\qquad\qquad \omega_3(\xi)\,>\,0.$$
Note the close analogy with the one dimensional case considered at (\ref{Taw}).
In turn, (\ref{Tlam}) yields
\bel{Ttau}\tau(\xi,\sigma)~=~-{1\over \lambda_1(\xi,\sigma)}~=~{1\over\omega_1(\xi)}
\left( 1 + {\omega_3(\xi)\over \omega_1(\xi)} \sigma^2+\O(1)\cdot \sigma^3\right).
\eeq

In addition, since all vectors $r_1(\xi,\sigma)$ have unit length and by definition $r_1(\xi,0)=r_1(\xi)$,
we have
\bel{r1r1}\bega{l}
\Big\langle r_1(\xi)~,~r_1(\xi,\sigma)-r_1(\xi)\Big\rangle~=~q_1(\xi)\sigma^2+q_2(\xi)\sigma^3+\O(1)\sigma^4,\\[4mm] \Big\langle r_1^\perp(\xi)~,~r_1(\xi,\sigma)-r_1(\xi)\Big\rangle~=~q_3(\xi)\,\sigma+q_4(\xi)\,\sigma^2+\O(1)\cdot \sigma^3,\enda\eeq
where
\bel{p1xi}q_1(\xi)=\Big\langle r_1(\xi)~,~ \frac12\frac{\partial^2r_1}{\partial\sigma^2}(\xi,0)\Big\rangle,\qquad
q_3(\xi)~\doteq~\Big\langle r_1^\perp(\xi)~,~{\partial r_1\over\partial \sigma}(\xi,0)\Big\rangle,
\eeq
and similar formulas hold for $q_2, q_4$.
We now compute
\bel{xxi}\bega{l}\ds x(\xi, Z)+ t  \bfw\bigl( x(\xi, Z)\bigr) - x(\xi,\sigma)- t\bfw\bigl(x(\xi,\sigma)\bigr)
~=~\int_\sigma^Z {d\over ds}
\Big[ x(\xi, s)+ t\, \bfw\bigl( x(\xi, s)\bigr) \Big]\, ds
\\[4mm]
\ds= ~\int_\sigma^Z    \bigl(1+ t \lambda_1(\xi,s) \bigr)r_1(\xi,s)\, ds
\\[4mm]
 \ds= ~\int_\sigma^Z
\bigl(t- \tau(\xi) \bigr) \lambda_1(\xi,s)r_1(\xi,s) \, ds +
\int_\sigma^Z     \tau(\xi)\Big[\lambda_1(\xi,s) -\lambda_1(\xi)\Big] r_1(\xi,s)\, ds
\,.\enda
\eeq
Taking the inner product with $r_1^\perp(\xi)$, using (\ref{r1r1}) and observing that
$$\lambda_1(\xi,s) - \lambda_1(\xi)~=~\frac{\omega_3(\xi)}{2}s^2+\frac{\omega_4(\xi)}{6}s^3+\frac{\omega_5(\xi)}{24}s^4+\O(1)\cdot s^5,$$ one obtains
\bel{P22}\bega{rl}
\Phi_2\bigl( t, \xi, \sigma, \xi, Z\bigr)&=\ds~\bigl( t-\tau(\xi)\bigr) \lambda_1(\xi) \Big[q_3(\xi)  {Z^2-\sigma^2\over 2}+
q_4(\xi)\frac{Z^3-\sigma^3}{3}\Big]+\frac{\omega_3(\xi)}{8}\tau(\xi)q_3(\xi)\bigl(Z^4-\sigma^4\bigr)\\[3mm]
& \qquad   \ds + \O(1) \cdot \bigl(t-\tau(\xi)\bigr)\bigl(Z^4-\sigma^4\bigr)+\,\O(1)\cdot\bigl(Z^5-\sigma^5\bigr).\enda
\eeq
Recalling (\ref{p1xi}) and (\ref{p2xi}), we set
\bel{q57}
 q_5(\xi)\doteq - {\lambda_1(\xi) \,q_3(\xi)\over 2 p_0(\xi) },\qquad q_6(\xi)\doteq-\frac{\lambda_1(\xi)q_4(\xi)}{3p_0(\xi)},
 \qquad q_7(\xi)\doteq-\frac{\omega_3(\xi)\tau(\xi)q_3(\xi)}{8p_0(\xi)}.\eeq
In view of (\ref{P22}) and (\ref{p2xi}),
by the implicit function theorem we conclude
\bel{phism}\bega{rl} \phi \bigl( t, \xi, \sigma, Z\bigr)-\xi &\ds=~
  \bigl( t-\tau(\xi)\bigr)\Big[ q_5(\xi)  (Z^2-\sigma^2)+q_6(\xi)(Z^3-\sigma^3)\Big]+q_7(\xi)\bigl(Z^4-\sigma^4\bigr)\\[4mm]
& \qquad   \ds + \O(1) \cdot \bigl(t-\tau(\xi)\bigr)\bigl(Z^4-\sigma^4\bigr)+\,\O(1)\cdot\bigl(Z^5-\sigma^5\bigr)\\[3mm]
&\doteq~\phi_0\bigl( t, \xi, \sigma, Z\bigr)~=~\O(1)\cdot \bigl(t-\tau(\xi)\bigr)^2\,.\enda
\eeq
Note that, in the case where the integral curves of $r_1$ are straight lines, that is: $r_1(\xi,\sigma)=r_1(\xi)$ for all $\sigma$,
we trivially have  $\phi_0(t,\xi,\sigma,Z)\equiv 0$. In particular,
$q_1(\xi)=\cdots =q_7(\xi)=0$.

Next,
inserting the value
 $$\xi_i~=~\phi(t,\xi, \sigma, Z) ~=~ \xi+ \phi_0(t,\xi,\sigma,Z)$$
 in the first equation of (\ref{le}), we recover a scalar
equation for the single
variable $Z\doteq \sigma_i$.  This has the form
\bel{eqZ}\bega{l}
\Tilde \Phi( t,\xi,\sigma,Z)~\doteq~\Phi_1\bigl(t,\xi,\sigma,\phi(t,\xi,\sigma,Z),Z
\bigr)\\[4mm]
\quad =~\ds\Big\langle r_1(\xi)\,,~ x\bigl(\phi(t,\xi,\sigma, Z), Z\bigr)+ t \, \bfw\bigl(  x\bigl(\phi(t,\xi,\sigma, Z), Z\bigr)\bigr) - x(\xi,\sigma)- t\, \bfw\bigl(x(\xi,\sigma)\bigr)
\Big\rangle~=~0.\enda\eeq
Notice that this is a scalar equation of cubic type, similar to (\ref{XXE}).
Indeed, when $\phi_0=0$, it reduces precisely to (\ref{XXE}).  Hence the solution
would be given by the same formula as (\ref{X2o}).

Using  the approximation (\ref{Tlam}) together with (\ref{r1r1}),  we
obtain
\bel{P1es}
\bega{l}\ds  \Tilde
\Phi(t,\xi,\sigma,\xi,Z)~\ds=~ \Big\langle r_1(\xi)\,,~ x(\xi, Z)+ t \, \bfw\bigl( x(\xi, Z)\bigr) - x(\xi,\sigma)- t\, \bfw\bigl(x(\xi,\sigma)\bigr)
\Big\rangle\\[4mm]
\qquad +
\Big\langle r_1(\xi)\,,~  x\bigl(\phi(t,\xi,\sigma, Z), Z\bigr)+ t \, \bfw\bigl(  x\bigl(\phi(t,\xi,\sigma, Z), Z\bigr)\bigr) - x(\xi, Z)- t \, \bfw\bigl( x(\xi, Z)\bigr)
\Big\rangle\\[4mm]
\quad = ~A+B.\enda\eeq
Recalling that $\omega_1(\xi) = -\lambda_1(\xi) = 1/\tau(\xi)$, we compute
\bel{A=}
\bega{rl}A&=~\ds \left\langle r_1(\xi)~,~
\int_\sigma^Z  \bigl(1+t  \lambda_1(\xi,s)\bigr) \, r_1(\xi,s)\, ds\right\rangle\\[4mm]
&=\ds\int_\sigma^Z \bigl(1+ t \lambda_1(\xi,s) \bigr)\, ds
\\[4mm]&\ds\qquad
+
 \left\langle r_1(\xi)~,~\int_\sigma^Z \bigl(1+t \lambda_1(\xi)\bigr)\bigl( r_1(\xi,s)-r_1(\xi)\bigr), ds
\right\rangle + \O(1)\cdot (Z-\sigma)\bigl(|Z|^4+|\sigma|^4\bigr)
\\[4mm]
&=~\ds -\omega_1(\xi) \bigl( t-\tau(\xi)\bigr) (Z-\sigma)+{t\,\omega_3(\xi)\over 6 }\,(Z^3-\sigma^3) +\frac{t\,\omega_4(\xi)}{24}\bigl(Z^4-\sigma^4\bigr)\\[4mm]
&\qquad \ds-\omega_1(\xi)  \bigl( t-\tau(\xi)\bigr) \biggl[\frac{q_1(\xi)}{3}\,(Z^3-\sigma^3)+
\frac{q_2(\xi)}{4}(Z^4-\sigma^4)\bigg]\\[4mm]
& \qquad +\O(1) \cdot \bigl(t-\tau(\xi)\bigr)\bigl(Z^4-\sigma^4\bigr)+\,\O(1)\cdot\bigl(Z^5-\sigma^5\bigr),
\enda
\eeq
where $q_1(\xi), q_2(\xi)$ are the coefficients in the Taylor expansion at (\ref{r1r1}).
Moreover, 
\bel{B=}
\bega{rl} B&=~\Big\langle r_1(\xi)\,,~\bigl(1+\tau(\xi) D\bfw(\xi)\bigr) r_2^*(\xi)
\Big\rangle \phi_0(t,\xi,\sigma,Z)
 \\[3mm]&\qquad +\O(1)\cdot \bigl(t-\tau(\xi)\bigr) \phi_0
 + \O(1)\cdot \phi_0^2 + \O(1)\cdot \phi_0 Z.
\enda
\eeq
To remove the trivial solution $Z=\sigma$,
we divide  (\ref{eqZ}) by $Z-\sigma$. In view of (\ref{phism}) we obtain
\bel{TP1}\bega{rl}\ds 0 &=~\ds
\Phi (t,\xi,\sigma,Z)~\doteq~{\Phi_1\bigl(t,\xi,\sigma,\phi(t,\xi,\sigma,Z),Z
\bigr)\over Z-\sigma}\\[4mm]
&=~\ds -\omega_1(\xi) \bigl( t-\tau(\xi)\bigr) +{t\,\omega_3(\xi)\over 6 }\,(Z^2+Z\sigma+\sigma^2) +\frac{t\,\omega_4(\xi)}{24}(Z+\sigma)(Z^2+\sigma^2)\\[4mm]
&\qquad \ds-\omega_1(\xi)  \bigl( t-\tau(\xi)\bigr) \biggl[\frac{q_1(\xi)}{3}\,(Z^2+Z\sigma+\sigma^2)+
\frac{q_2(\xi)}{4}(Z+\sigma)(Z^2+\sigma^2)\bigg]\\[4mm]
& \qquad + q_8(\xi)\Big\{\bigl( t-\tau(\xi)\bigr)\Big[ q_5(\xi)(Z+\sigma)+q_6(\xi)(Z^2+\sigma Z+\sigma^2)\Big]
+q_7(\xi)(Z+\sigma)\bigl(Z^2+\sigma^2\bigr)\Bigr\}
\\[4mm]&\ds
\qquad\qquad + \O(1) \cdot \bigl(t-\tau(\xi)\bigr)\bigl(|Z|^3+|\sigma|^3\bigr)+\,\O(1)\cdot\bigl(|Z|^4+|\sigma|^4\bigr),\enda
 \eeq
 where
 $$q_8(\xi)~\doteq~\Big\langle r_1(\xi)\,,~\bigl(I+\tau(\xi)D_x\bfw(\xi)\bigr)r_2^*(\xi)
\Big\rangle .$$

Dividing all terms by $t\omega_3(\xi)/6$ and recalling that $\tau(\xi) \omega_1(\xi)=1$,
we are led to the equation
\bel{EZ1}\bega{l}\ds 0~=~
f(t,\xi,\sigma,Z)~\doteq~(Z^2+Z\sigma+\sigma^2)
\ds -{6 \omega_1(\xi) \bigl( t-\tau(\xi)\bigr)\over t\, \omega_3(\xi)}
+\frac{\omega_4(\xi)}{4\omega_3(\xi)}(Z+\sigma)(Z^2+\sigma^2)\\[4mm]
\qquad \ds-{6\omega_1(\xi)  \bigl( t-\tau(\xi)\bigr)\over t\,\omega_3(\xi)} \biggl[\frac{q_1(\xi)}{3}\,(Z^2+Z\sigma+\sigma^2)+
\frac{q_2(\xi)}{4}(Z+\sigma)(Z^2+\sigma^2)\bigg]\\[4mm]
 \qquad \ds+ {6q_8(\xi)\over t\,\omega_3(\xi)}  \left\{ \bigl( t-\tau(\xi)\bigr)\Big[ q_5(\xi)(Z+\sigma)+q_6(\xi)(Z^2+\sigma Z+\sigma^2)\Big]
+q_7(\xi)(Z+\sigma)\bigl(Z^2+\sigma^2\bigr)\right\}
\\[4mm]\ds
\qquad+ \E(t,\xi,\sigma,Z),\enda
 \eeq
 where $\E$ is an analytic function that collects higher order terms:
 \bel{gho} \E(t,\xi,\sigma,Z)~\doteq~
 \O(1) \cdot \bigl(t-\tau(\xi)\bigr)\bigl(|Z|^3+|\sigma|^3\bigr)+\,\O(1)\cdot\bigl(|Z|^4+|\sigma|^4\bigr).
 \eeq
Note that, if the integral curves of the eigenvector $r_1$ are straight lines, i.e.
$r_1(\xi,\sigma)= r_1(\xi)$, then $q_1(\xi)=\cdots=q_8(\xi)=0$ and  above formula reduces to (\ref{xe1}), with $t_0=\tau(\xi)$.

Since  (\ref{EZ1}) will be solved on a domain of the form
\bel{Dom} \D~=~\Big\{ (t,\xi,\sigma)\,;~~0<t-\tau(\xi)<\ve,~~|\sigma|\leq M\, \bigl(t-\tau(\xi)\bigr)\Big\},\eeq
it will be useful to make the change of variable
\bel{Tsig}
\Tilde\sigma~=~{\sigma\over t-\tau(\xi)}\,.
\eeq
Still in view of (\ref{Dom}), we expect $Z= \O(1)\cdot \bigl(t-\tau(\xi)\bigr)^{1/2}$.
Using the elementary expansion
$${1\over t}~=~{1\over\tau(\xi) + \bigl(t-\tau(\xi)\bigr)}~=~
{1\over\tau(\xi)} \cdot \sum_{k=0}^\infty \left(t-\tau(\xi)\over\tau(\xi)\right)^k,\qquad\qquad {1\over\tau(\xi)}~=~\omega_1(\xi),$$
and arranging terms according to their order of magnitude,  we write (\ref{EZ1}) in the form
\bel{Z3}\bega{rl}\ds 0&\ds=~\bigg\{\Big[
Z^2 - c_{1}(\xi) \bigl(t-\tau(\xi)\bigr) \Big]  +  \Big[ c_2(\xi,\Tilde\sigma)
 \bigl(t-\tau(\xi)\bigr) Z +c_3(\xi)
 Z^3  \Big]  \bigg\}
 \\[4mm]&\qquad \qquad \ds
 +\sum_{j+(k/2)\geq 2,\;j>0}
c_{jk}(\xi,\Hat \sigma)  \bigl(t-\tau(\xi)\bigr)^j Z^k  \\[4mm]
\qquad &\doteq~\Gamma(t,\xi,\sigma,Z)+\E(t,\xi,\sigma,Z).
\enda
\eeq
Here
\bel{bb}c_{1}(\xi)~=~ {6 \omega_1^2(\xi) \over \omega_3(\xi)}\,,\qquad\qquad
 c_2(\xi,\Tilde\sigma)\,=\,\Tilde\sigma+{6  q_5(\xi)q_8(\xi)\omega_1(\xi) \over  \omega_3(\xi)} ,\eeq
\bel{cc} c_3(\xi)\,=\,\frac{\omega_4(\xi)}{4\omega_3(\xi)}+\frac{6q_7(\xi)q_8(\xi)\omega_1(\xi)}{\omega_3(\xi)}\,,\eeq
while $c_{jk}(\xi,\Tilde \sigma)$ are analytic functions.

\begin{lemma}\label{l:41} Given $\xi_0, \Tilde\sigma_0$,
there exist analytic functions $\beta_0, \beta_1$ defined in a neighborhood of
$\bigl( \tau(\xi_0), \xi_0, \Tilde\sigma_0\bigr)$ such that the two solutions to
(\ref{Z3}) have the form
\bel{Zrep}
Z^\pm(t,\xi,\Tilde\sigma)~=~\bigl(t-\tau(\xi)\bigr)\cdot\beta_0(t,\xi,\Tilde\sigma)\, \pm\, \sqrt{ t-\tau(\xi)} \cdot \beta_1(t,\xi,\Tilde\sigma).\eeq
The above coefficients satisfy
\bel{bet0}
\beta_0(t,\xi,\Tilde\sigma)~=~-{c_2(\xi,\Tilde\sigma)+ c_1(\xi)c_3(\xi)\over 2} +
\O(1)\cdot \bigl(t-\tau(\xi)\bigr) ,\eeq
\bel{bet1}\beta_1(t,\xi,\Tilde\sigma)~=~\sqrt{c_1(\xi)}+\O(1)\cdot\bigl(t-\tau(\xi)\bigr).\eeq
\end{lemma}

\begin{remark} {\rm
It is worth noting that the formulas (\ref{Zrep})--(\ref{bet1}) reduce to  (\ref{X2o})-(\ref{BCdef}) in the case where
the integral curves of the eigenvectors $r_1$ are straight lines,  so that
$q_1(\xi)=\cdots=q_8(\xi)=0$.
%
}
\end{remark}

{\bf Proof of Lemma~\ref{l:41}.}

 {\bf 1.} To construct a solution to (\ref{Z3}), we
start with the approximation
$$Z~\approx~\pm\sqrt{c_1(\xi) \bigl(t-\tau(\xi)\bigr)}.$$
To the next order, inserting the above value to compute $Z^3$,
we obtain an approximate solution to
$\Gamma(t,\xi,\sigma,Z)=0$, namely
\bel{Zpm0}\bega{rl}
Z_0^\pm(t,\xi,\Tilde\sigma)&\ds\doteq \,-{\bigl(t-\tau(\xi)\bigr)c_2(\xi,\Tilde\sigma) \over 2}\\[4mm]
&\qquad\ds
 \pm {1\over 2}
\sqrt{ \bigl[\bigl(t-\tau(\xi)\bigr)c_2(\xi,\Tilde\sigma) \bigr]^2 +4 c_1(\xi) \bigl(t-\tau(\xi)\bigr)\mp
4c_3(\xi)\bigl[ c_1(\xi)(t-\tau(\xi))\bigr]^{3/2}}.\enda
\eeq
We observe that the argument of the square root is a small perturbation of its leading term
$4 c_1(\xi) \bigl(t-\tau(\xi)\bigr)$.  Expanding this argument in Taylor series, and
collecting even and odd powers of $\bigl(t-\tau(\xi)\bigr)^{1/2}$, one can write
(\ref{Zpm0}) in the form
\bel{Zpm1}
Z_0^\pm(t,\xi,\Tilde\sigma)~=~b(t,\xi,\Tilde\sigma)\cdot \bigl(t-\tau(\xi)\bigr)\pm h(t,\xi,\Tilde\sigma) \cdot \sqrt{t-\tau(\xi)},\eeq
for suitable functions $b,h$, analytic on a neighborhood of $\bigl( \tau(\xi_0), \xi_0, \Tilde\sigma_0\bigl)$.
A more detailed expansion yields
\bel{Zpm2}\bega{rl}
Z_0^\pm(t,\xi,\Tilde\sigma)
&\ds=~
-{\bigl(t-\tau(\xi)\bigr)c_2(\xi,\Tilde\sigma) \over 2}\pm
\sqrt{c_1(\xi)\bigl(t-\tau(\xi)\bigr)}\cdot\sqrt{1\mp c_3(\xi)
\left[c_1(\xi)\bigl(t-\tau(\xi)\bigr)\right]^{1/2}}\\[4mm]
&\qquad
+\O(1)\cdot\bigl(t-\tau(\xi)\bigr)^{3/2}\\[4mm]
&\ds=
\,-{\bigl(t-\tau(\xi)\bigr)c_2(\xi,\Tilde\sigma) \over 2}\pm
\sqrt{c_1(\xi)\bigl(t-\tau(\xi)\bigr)}- \frac{c_1(\xi)}{2} c_3(\xi)
\bigl(t-\tau(\xi)\bigr)\\[4mm]
&\qquad+\O(1)\cdot\bigl(t-\tau(\xi)\bigr)^{3/2}.
\enda
\eeq
Therefore, in (\ref{Zpm1}) one has
\bel{be0}b(t,\xi,\Tilde\sigma)~=~-{c_2(\xi,\Tilde\sigma)+ c_1(\xi)c_3(\xi)\over 2}
 +
\O(1)\cdot \bigl(t-\tau(\xi)\bigr) ,\eeq
\bel{be1}h(t,\xi,\Tilde\sigma)~=~\sqrt{c_1(\xi)}+\O(1)\cdot\bigl(t-\tau(\xi)\bigr).\eeq
Recalling (\ref{bb})-(\ref{cc}), we see that the above formulas reduce to (\ref{X2o})-(\ref{BCdef}) if the integral curves of $r_1(x)$ are straight lines, so that $q_1(\xi)=\cdots=q_8(\xi)$.

\v
{\bf 2.} Starting with $Z_0$ as a first approximation,
and applying Newton's iterative algorithm,
we obtain a sequence of approximations inductively defined as
\bel{Newton}
Z_{n+1}~=~Z_n - {(\Gamma+\E)(t,\xi,\Tilde\sigma, Z_n)\over \partial_Z (\Gamma+\E)(t,\xi,\Tilde\sigma, Z_n)}\,.\eeq
We claim that
this sequence is uniformly convergent on a neighborhood of $\bigl( \tau(\xi_0), \xi_0, \Tilde\sigma_0\bigl)$.
Indeed
\begi
\item[(i)] 
By construction,
the functions $Z_0^\pm$ satisfy an equation of the form
\bel{co1}\bega{rl}
(\Gamma+\E)(t,\xi, \Tilde \sigma, Z_0^\pm)&\ds = \sum_{j+(k/2)\geq 2,\, j>0}
c_{jk}(\xi,\Tilde\sigma)  \bigl(t-\tau(\xi)\bigr)^j (Z_0^\pm)^k ~=~\O(1)\cdot  \bigl(t-\tau(\xi)\bigr)^2 .\enda
 \eeq
\item[(ii)] At the  point $(t,\xi,  \Tilde \sigma,
Z_0^\pm)$ the partial derivative $ \partial_Z(\Gamma+\E)$
is bounded away from zero.
Indeed, differentiating (\ref{Z3}) w.r.t.~$Z$ one obtains
$$\partial_Z(\Gamma+\E)(t,\xi,\Tilde \sigma, Z)~=~2Z+
\sum_{j+(k/2)\geq 2}
c_{jk}(\xi,\Tilde\sigma)  \bigl(t-\tau(\xi)\bigr)^j\cdot k \, Z^{k-1}.$$
Therefore, for $(t,\xi,\Tilde \sigma)$ in a neighborhood of  $\bigl( \tau(\xi_0), \xi_0, \Tilde\sigma_0\bigl)$,
 there exists a constant $\kappa_0>0$ such that
\bel{co2}\Big|\partial_Z(\Gamma+\E)(t,\xi,\Tilde\sigma, Z_0^\pm)\Big|
~\geq~\kappa_0 \bigl(t-\tau(\xi)\bigr)^{1/2}.\eeq

\item[(iii)] Finally, from (\ref{Z3}) it is clear that the second derivative $\partial_{ZZ} (\Gamma+\E)(t,\xi,\Tilde \sigma, Z)$
is locally bounded
\endi
By Kantorovich' theorem (see for example \cite{Berger, Deimling}),
the sequence $(Z_n^\pm)_{n\geq 0}$ converges to a limit
$Z^\pm=Z^\pm(t,\xi,\Tilde\sigma)$ which provides a solution to (\ref{Z3}).
Moreover
\bel{ZZZ}
\bigl| Z^\pm - Z_0^\pm\bigr| ~=~\O(1)\cdot {(\Gamma+\E)(t,\xi,\Tilde\sigma, Z_0^\pm)\over \partial_Z (\Gamma+\E)(t,\xi,\Tilde\sigma, Z_0^\pm)}
~\leq~\kappa_1 \cdot \bigl(t-\tau(\xi)\bigr)^{3/2},
\eeq
 for some constant $\kappa_1$, uniformly valid on a neighborhood of  $\bigl( \tau(\xi_0), \xi_0, \Tilde\sigma_0\bigl)$.
\v
{\bf 3.} It remains to check the regularity of the maps $Z^\pm (t,\xi,\Tilde\sigma)$.
Since the functions $\Gamma$, $\E$ are analytic, every iterate $Z_n^\pm$ can be written
as a power series in $Z_0^\pm$,
$$Z_n^\pm~=~\sum_j c_j^{(n)} (t,\xi,\Tilde\sigma) (Z_0^\pm)^j.$$
By uniform convergence, the same holds for the limit:
$$Z^\pm~=~\sum_j c_j(t,\xi,\Tilde\sigma) (Z_0^\pm)^j.$$
Recalling (\ref{Zpm1}), we can separate terms containing  even or odd powers of
$\sqrt{t-\tau(\xi)}$.   This leads to (\ref{Zrep}).
\v
{\bf 4.}
Finally, the identities in (\ref{bet0}-(\ref{bet1})
follow from (\ref{be0})-(\ref{be1}).
\endproof
\v
In turn,  (\ref{Zrep}) yields information about the velocity $\bfv^\pm$
and density function $\rho^\pm$, in a neighborhood of a point where the singularity is formed.
As in (\ref{defs}), in the following we continue to use the shorter notation $\lambda_1(\xi)\doteq\lambda_1\bigl(\gamma_1^*(\xi)\bigr)$, ~
$\bar\rho(\xi)\doteq\bar\rho\bigl(\gamma_1^*(\xi)\bigr)$, etc$\ldots$

\begin{lemma}\label{l:42} Let $\xi_0$ and $\Tilde\sigma_0\in\R$ be given.   Then there exist
analytic functions  $\bfB_1, \bfB_2, b_1, b_2$, defined for
$(t,\xi, \Tilde\sigma)$ in a neighborhood of $\bigl( \tau(\xi_0), \xi_0, \Tilde\sigma_0\bigr)$,
such that the following  holds.

 For
$t\geq \tau(\xi)$ and $\sigma = \bigl(t-\tau(\xi)\bigr)\Tilde\sigma$, one has
\bel{vpm2}
\bfv^\pm(t,\xi,\sigma)-\bfw(\xi)~=~\bfB_1(t,\xi,\Tilde\sigma) \cdot  \bigl(t-\tau(\xi)\bigr)\pm \bfB_2(t,\xi,\Tilde\sigma)\cdot \sqrt{t-\tau(\xi)}\,,
\eeq
\bel{rpm2}
\rho^{\pm}(t,\xi,\sigma)~ =~\frac{1}{t-\tau(\xi)}\Big[b_1(t,\xi,\Tilde\sigma) \pm b_2(t,\xi,\Tilde\sigma)\cdot \sqrt{t-\tau(\xi)}\Big].
 \eeq
Moreover, recalling (\ref{bb})-(\ref{cc}) one obtains
 \bel{bB0} 
 \bfB_2(\tau(\xi),\xi,\Tilde\sigma)\,=\,\sqrt{c_1(\xi)}\lambda_1(\xi)r_1(\xi)~=~-\omega_1^2(\xi) \sqrt{6\over \omega_3(\xi) } \, r_1(\xi)\,\not=\,0,
\eeq
 \bel{b12}
 b_1(\tau(\xi),\xi,\Tilde\sigma)~=~\frac{\bar\rho(\xi)}{2(\lambda_2(\xi)-\lambda_1(\xi))}\,>\,0
.\eeq
 \end{lemma}

{\bf Proof.} {\bf 1.} The representation (\ref{vpm2}) follows from
\bel{v22}\bfv^\pm(t,\xi,\sigma)~=~\bfw\bigl(X^\pm(t,\xi,\sigma)\bigr),\eeq
where
\bel{v33}X^\pm(t,\xi,\sigma)~=~x\Big( \phi(t,\xi,\sigma, Z^\pm), \, Z^\pm\Big)
\eeq
are the initial points of the two additional
characteristics reaching the point $x(\xi,\sigma) + t \bfw\bigl( x(\xi,\sigma)\bigr)$
at time $t$.  

In view of Lemma~\ref{l:41}, the two functions $Z^\pm(t,\xi, \Tilde\sigma)$ admit the analytic representation
(\ref{Zrep}).
Since the map $(\xi,\sigma)\mapsto x(\xi,\sigma) $ at (\ref{xxs})
is analytic, and so is the map $\phi$ at (\ref{xii}), it is clear that the two
composed maps in (\ref{v22}) can be represented as in (\ref{vpm2}).

Recalling (\ref{phism}) and (\ref{bet1}), when $t=\tau(\xi)$ and hence
$\sigma = \bigl(t-\tau(\xi)\bigr) \Tilde\sigma=0$,
we obtain
$$\bfB_2\bigl(\tau(\xi),\xi,\Tilde\sigma \bigr)~=~\beta_1\bigl(\tau(\xi),\xi,\Tilde\sigma \bigr) \,D \bfw(\xi) r_1(\xi)
~=~\sqrt{c_1(\xi)}\lambda_1(\xi)r_1(\xi)~\not=~0.$$
In addition, since $\bfB_1$ collects even powers of $\sqrt{t-\tau(\xi)}$ and $\bfB_1
\bigl(\tau(\xi),\xi,\Tilde\sigma\bigr)=0$,
by analyticity
the second estimate in (\ref{bB0}) must hold.

Furthermore, by Taylor expansions we have the following leading order approximations:
\bel{v44}\bega{l}\bfB_1(t,\xi,\Tilde \sigma)~\ds=~
\beta_0(t,\xi,\Tilde\sigma)\lambda_1(\xi)r_1(\xi)\\[4mm]
\qquad\ds +\frac{c_1(\xi)}{2}\Bigl(D^2\bfw(\xi) \bigl(r_1(\xi)\otimes r_1(\xi)\bigr)+D\bfw(\xi)\,\bigl(Dr_1(\xi)\cdot r_1(\xi)\bigr)\Bigr)
 +\O(1)\cdot\bigl(t-\tau(\xi)\bigr),\enda\eeq
\bel{v55}\bfB_2(t,\xi,\Tilde \sigma)=\sqrt{c_1(\xi)}\lambda_1(\xi)r_1(\xi)
+\O(1)\cdot\bigl(t-\tau(\xi)\bigr).\eeq
\v
{\bf 2.}
To prove  (\ref{rpm2}) we recall that, by (\ref{den}) and
(\ref{det12}),
\bel{rhopp}
\rho^\pm (t,\xi,\sigma)~=~{\bar \rho(X^\pm)\over \det\bigl(I + tD\bfw(X^\pm)\bigr)}~
=~{\bar \rho(X^\pm)\over
1+t\lambda_2\bigl(X^\pm\bigr) }
\cdot {1\over 1+t\lambda_1\bigl(X^\pm)}
\,,\eeq
with $X^\pm$ as in (\ref{v33}).
The same arguments as in step {\bf1} yield the existence of scalar functions
$\vp_1, \vp_2, \vp_3,\vp_4$, analytic in a neighborhood of  $\bigl( \tau(\xi_0), \xi_0, \Tilde\sigma_0\bigr)$,
such that
\bel{te1}{\bar \rho(X^\pm)\over
1+t\lambda_2\bigl(X^\pm\bigr) }~=~\vp_1(t,\xi,\Tilde\sigma)\pm \sqrt{ t-\tau(\xi)}\cdot \vp_2(t,\xi,\Tilde\sigma),
\eeq
\bel{te2}
 1+t\lambda_1\bigl(X^\pm)~=~\bigl(t-\tau(\xi)\bigr)\cdot\vp_3(t,\xi,\Tilde\sigma)\pm \sqrt{ t-\tau(\xi)} \cdot\vp_4(t,\xi,\Tilde\sigma).\eeq
 To estimate the size of (\ref{te2}), and understand at what rate its inverse blows up,
 we compute
\bel{te3} 1+t\lambda_1\bigl(X^\pm)~=~-\lambda_1(\xi) \tau(\xi) +t\lambda_1\bigl(X^\pm)
~=~ \bigl(t-\tau(\xi)\bigr) \lambda_1(\xi)+ t \Big[ \lambda_1(X^\pm) -\lambda_1(\xi)\Big].
\eeq
Recalling (\ref{Tlam}),
 (\ref{v33}), and the last estimate
in (\ref{phism}),
 we obtain
\bel{te4}  \lambda_1(X^\pm) -\lambda_1(\xi)~=~{\omega_3(\xi)\over 2}( Z^\pm)^2 + \O(1) \cdot  \bigl(t-\tau(\xi)\bigr)^{3/2}.
\eeq

Using  Lemma~\ref{l:41}  to compute $Z^\pm$, and recalling (\ref{bb}) and the identity
$\tau(\xi) \omega_1(\xi)=1$, one obtains
\bel{55}\bega{l} \ds
\bigl(t-\tau(\xi)\bigr) \lambda_1(\xi)+ t \Big[ \lambda_1(X^\pm) -\lambda_1(\xi)\Big]\\[4mm]
\qquad\ds
=~  -\omega_1(\xi) \bigl(t-\tau(\xi)\bigr) +\tau(\xi)\cdot {\omega_3(\xi)\over 2}\,c_1(\xi) \bigl(t-\tau(\xi)\bigr)+ \O(1)\cdot (t-\tau(\xi))^{3/2}\\[4mm]
\qquad\ds =~\left[ -\omega_1(\xi) + {1\over\omega_1(\xi)} \cdot {\omega_3\over 2}
\cdot {6\omega_1^2(\xi)\over \omega_3(\xi)}\right]\bigl(t-\tau(\xi)\bigr)+ \O(1)\cdot \bigl(t-\tau(\xi)\bigr) ^{3/2}\\[4mm]
=~\ds 2 \omega_1(\xi) \bigl(t-\tau(\xi)\bigr)+ \O(1)\cdot \bigl(t-\tau(\xi)\bigr) ^{3/2} .\enda
\eeq
Therefore, the inverse of the quantity in (\ref{te3}) can be written as
\bel{66}
{1\over 1+t\lambda_1\bigl(X^\pm)}~=~{1\over 2\omega_1(\xi) \bigl(t-\tau(\xi)\bigr) }\cdot \left(
1+\sum_{k\geq 1} d_k \bigl(t-\tau(\xi)\bigr)^{k/2}\right).\eeq
Together with (\ref{te1}), this yields
(\ref{rpm2}).
Since we  assume that the initial density
$\bar \rho$ is uniformly positive,
taking the limit as $t\to \tau(\xi)+$, the previous analysis yields
\bel{b1e}b_1\bigl( \tau(\xi), \xi,\Tilde\sigma\bigr)~=~\frac{\tau(\xi)}{2}\cdot{\bar\rho(\xi)\over 1+ \tau(\xi)\lambda_2(\xi)}~=~\frac{\bar\rho(\xi)}{2\bigl(\lambda_2(\xi)-\lambda_1(\xi)\bigr)}
~>~0.\eeq
Indeed, $\tau(\xi)=-1/\lambda_1(\xi)$ and $\lambda_2(\xi)>\lambda_1(\xi)$.
\endproof

\section{Local construction of the singular curve}
\label{s:5}
\setcounter{equation}{0}

In view of (\ref{92}), (\ref{94}), we
seek a solution $(y,\eta, \bfv)(t,\zeta)\in \R^2\times\R\times\R^2$ of the system
\bel{2ds}
\left\{\bega{rl}\eta_t&=~
 \rho^+ (\bfv-\bfv^+)\times y_\zeta +\rho^- (\bfv^--\bfv)\times y_\zeta\,,\\[2mm]
y_t&=~\bfv,\\[2mm]
\bfv_t&=~\ds
\rho^+\Big( (\bfv-\bfv^+)\times y_\zeta\Big)\, {\bfv^+-\bfv\over\eta} +  \rho^-\Big( (\bfv^--\bfv)\times y_\zeta\Big)\, {\bfv^--\bfv\over\eta}  \,.\enda\right.\eeq
The boundary data are
\bel{bda1}
\left\{\bega{rl} y(\tau(\zeta),\zeta) & =~\gamma^\sharp(\zeta)~\doteq~\gamma^*(\zeta)+ \tau(\zeta) \bfw(\gamma^*(\zeta)),\\[2mm]
\eta(\tau(\zeta),\zeta) & =~0,\\[2mm]
\bfv(\tau(\zeta),\zeta)&=~\bfw(\gamma^*(\zeta)).\enda\right.\eeq
Here we use  $t,\zeta$ as independent variables to describe
points on the singular curve (see Fig.~\ref{f:gen60}). For a fixed $\zeta$,
the map  $t\mapsto y(t,\zeta)$, with  $t\geq \tau(\zeta)$,
describes the trajectory of a point on the singular
curve.  More precisely,
consider the  particle that  at time $t=0$ was located at $\gamma^*(\zeta)$.
\begi
\item During the time interval $t\in \bigl[0, \tau(\zeta)\bigr] $, this particle moves with constant speed $\bfw\bigl(\gamma^*(\zeta)\bigr)$.
\item At time  $t=\tau(\zeta)$, this particle is located at $\gamma^\sharp(\zeta)$, which is
one of the endpoints of the singular curve $\gamma^t$.
\item At times
 $t>\tau(\zeta)$ this same particle is located at a point $y(t,\zeta)$ on the singular
 curve~$\gamma^t$. \endi
Notice that, at a fixed time $t>t_0$, the singular curve is parameterized
by the map $\zeta\mapsto y(t,\zeta)$. Here $\zeta$ ranges over the set
$\bigl\{ \zeta\,;~\tau(\zeta)\leq t\bigr\}$.

To state the main result, providing an asymptotic description of the
initial stages of the singular curve, we need to work in
a suitable space of analytic vector fields.
More precisely, we fix a closed ball ${\cal B}= \{ x\in \R^2\,;~|x|\leq r_0\} $
and a constant
$\varrho>0$.   We call $X_{{\cal B}, \varrho}$ the space of all
analytic functions $\phi:\R^2\mapsto\R$ such that
\bel{avf}   \|\phi\|_{{\cal B}, \varrho}~\doteq~\sum_{n=0}^\infty {\varrho^n\over n!}
\max_{|\alpha|=n} \sup_{x\in {\cal B}} \bigl| D^\alpha \phi(x)\bigr|~<~+\infty.\eeq
Here $\alpha=(\alpha_1,\alpha_2)$ is a multi-index, while
$D^\alpha \phi\doteq \partial_{x_1}^{\alpha_1}\partial_{x_2}^{\alpha_2}\phi$.
Notice that this is the set of all functions whose Taylor series at every point $x\in {\cal B}$ converges on a disc
of radius $\varrho/2$.

In the following theorem, we consider the space of all initial data
$(\bfw,\bar\rho)= (w_1, w_2,\bar \rho)$ such that $w_1, w_2,\bar \rho\in X_{{\cal B}, \varrho}$, for some given ${\cal B}, \varrho$.   Within this space, we will construct
piecewise analytic solutions to the equations of pressureless gases (\ref{SP})-(\ref{ic})
beyond singularity formation, for a generic set of initial data.

\begin{theorem}\label{t:51} For the Cauchy problem (\ref{SP})-(\ref{ic}), assume that the initial data
are analytic functions, satisfying {\bf (A1)-(A2)}.

Then an  analytic solution up to some time $t_0$ where the $\C^1$ norm of the velocity $\bfv$ blows up.
This solution can be extended to a larger interval  $[0, t_0+\ve_0]$. For $t_0<t\leq t_0+\ve_0$
the solution is piecewise analytic, and contains a singular curve $\gamma^t$ where a positive amount of mass is concentrated.
\end{theorem}

{\bf Proof.} The analysis of solutions up to the first time $t_0$ when a singularity
appears, so that $\bigl\|\bfv(t,\cdot)\bigr\|_{\C^1}\to +\infty$ and $\bigl\|\rho(t,\cdot)\bigr\|_{\L^\infty}\to +\infty$, has been worked out in the recent paper \cite{BCH}.
The following proof will focus on the construction of a local
solution for $t\geq t_0$.
To help the reader, we summarize the main steps.
\begi
\item[(i)]
By a suitable change of variables,  the equations in (\ref{2ds})-(\ref{bda1})
are written as a system of 5 first order equations
(\ref{etat})--(\ref{Hsi}), with initial data (\ref{IC0})-(\ref{IC1}).
\item[(ii)] In turn, this system can be put in the standard form (\ref{sip}).  An integrability condition
near $T=0$ allows us to uniquely determine the initial  values  $\rho_0, S_0, Z_0$ in (\ref{IC1}).
\item[(iii)] By checking that none of the eigenvalues of the $5\times 5$ matrix $ {\bf\Lambda}(\zeta)=\frac{\partial \bfH}{\partial \bfu}(0,\zeta,0,0)$ in (\ref{chm}) is a positive integer,
we show that the system (\ref{sip}) is of Briot-Bouquet type.  The general result in \cite{li} can thus be applied, yielding the local existence of an analytic solution.
\endi

\v
{\bf 1.} As in (\ref{yxs}), a point $(t,y)\in [0,T]\times \R^2$ will be described
in terms of  the variables $(t,\xi,\sigma)$.  These are uniquely determined by the requirements
\bel{ytz} y~=~
x(\xi,\sigma) + t \bfw\bigl(x(\xi,\sigma)\bigr),
\qquad\qquad\tau\bigl(x(\xi,\sigma)\bigr) < t.\eeq
The main advantage of these variables is that now we can use the expressions
(\ref{v33}) for the initial points of the additional two characteristics that meet at $(t,y)$.
In turn, the formulas (\ref{vpm2})-(\ref{rpm2}) in Lemma~\ref{l:42} yield an accurate asymptotic
 description of  the
functions $\bfv^\pm,\rho^\pm$, near the singularity.
\v
Starting from the system (\ref{2ds}), we seek a system of equations for the
variables $(\xi, \sigma, \eta)$ as functions of $(t,\zeta)$, where $(t, x(\xi,\sigma))\in R$ is
related to $(t,y)$ by (\ref{yxs}).
Differentiating (\ref{ytz}), one obtains
\bel{yte}\left\{\bega{rl}
y_t&=~\bfw+ (I+t D_x\bfw )\cdot  (x_\xi \xi_t + x_\sigma \sigma_t),\\[3mm]
y_\zeta&=~ (I+t D_x\bfw )\cdot  (x_\xi \xi_\zeta + x_\sigma \sigma_\zeta),\\[3mm]
y_{tt}&=~2D_x\bfw \cdot  (x_\xi \xi_t + x_\sigma \sigma_t) +
 (x_\xi \xi_t + x_\sigma \sigma_t)^T(t\, D^2\bfw ) (x_\xi \xi_t + x_\sigma \sigma_t)\\[4mm]
&\qquad +  (I+t D_x\bfw )\cdot  \bigl(x_\xi \xi_{tt} + x_\sigma \sigma_{tt}+
x_{\xi\xi} \xi_t^2 + 2 x_{\xi\sigma} \xi_t \sigma_t + x_{\sigma\sigma} \sigma_t^2\bigr).
\enda\right.
\eeq
%
The first equation in (\ref{2ds}) yields
{\small
\bel{ett}
\bega{rl} \eta_t
&=~\bigg\{
\rho^+(\bfw-\bfv^+) + \rho^-(\bfv^--\bfw)
+
 (I+t D_x\bfw )\cdot  (x_\xi \xi_t + x_\sigma \sigma_t)(\rho^+ -\rho^-)\bigg\}\\[3mm]
 &\qquad\qquad \times (I+t D_x\bfw )\cdot  (x_\xi \xi_\zeta + x_\sigma \sigma_\zeta)\\[3mm]
&\doteq~H ( t,\xi, \sigma, \xi_t,\sigma_t, \xi_\zeta, \sigma_\zeta).
\enda
\eeq
}
Since $\bfv_t = y_{tt}$, by (\ref{yte}), the last equation in (\ref{2ds}) yields
\bel{ytt}
\bega{l}  (I+t D_x\bfw )\cdot  (x_\xi \xi_{tt} + x_\sigma \sigma_{tt})\\[4mm]
\qquad =~- 2D_x\bfw \cdot  (x_\xi \xi_t + x_\sigma \sigma_t) -
 (t\, D^2\bfw )\bigl( (x_\xi \xi_t + x_\sigma \sigma_t)\otimes(x_\xi \xi_t + x_\sigma \sigma_t)\bigr) \\[4mm]
\qquad \qquad -  (I+t D_x\bfw )\cdot  (
x_{\xi\xi} \xi_t^2 + 2 x_{\xi\sigma} \xi_t \sigma_t + x_{\sigma\sigma} \sigma_t^2)\\[4mm]
\qquad\qquad
\ds + {1\over\eta}\bigg\{
\bigl( (y_t-\bfv^+)\times y_\zeta \bigr) \rho^+ (\bfv^+ - y_t)  +
\bigl( (\bfv^--y_t) \times y_\zeta\bigr) \rho^- (\bfv^- - y_t)\bigg\}
 \\[4mm]
\ds\qquad \doteq~F( t, \xi, \sigma, \xi_t, \sigma_t)+ {1\over\eta} G( t,\xi, \sigma, \xi_t,
\sigma_t, \xi_\zeta, \sigma_\zeta),\enda\eeq
where $y_\zeta, y_t$ are also computed by the first two equations in (\ref{yte}).

\v
 {\bf 2.} In order to replace the vector equation (\ref{ytt}) with a couple of scalar equations, we
 consider the vectors
\bel{r12*}r_1(\xi,\sigma)~\doteq~{\partial\over \partial \sigma} x(\xi,\sigma),\qquad\qquad
 r_2^*(\xi,\sigma)~\doteq~{\partial\over \partial \xi} x(\xi,\sigma).\eeq
 By (\ref{xxs}), $r_1(\sigma,\xi)$ is a unit  eigenvector of the Jacobian matrix $D\bfw$ at the point $x(\xi,\sigma)$, corresponding to the first eigenvalue $\lambda_1$.
 However, $r_2^*$ is not an eigenvector, in general.      We denote by $r_2$ the unit eigenvector of $D\bfw$ corresponding to the second eigenvalue $\lambda_2$.
  In addition, we denote by
$\{ l_1, l_2\} $ the dual basis of $\{r_1, r_2\}$, and by
$\{ l_1^*, l_2^*\} $ the dual basis of $\{r_1, r_2^*\}$ so that
$$\langle l_1^*, x_\sigma\rangle \,=\,\langle l_2^*, x_\xi\rangle\,=\,1,  \qquad\qquad
 \langle l_1^*, x_\xi\rangle \,=\,\langle l_2^*, x_\sigma\rangle\,=\,0.$$

 Notice that $l_1,l_2$ are left eigenvectors of $D_x\bfw$ corresponding to the eigenvalues $\lambda_1< \lambda_2$, respectively.  Moreover, $l_2^*$ is also
 a left eigenvector with eigenvalue $\lambda_2$.
 However,  $l_1^*$ is not an eigenvector,
 in general.  Multiplying both sides of (\ref{ytt}) by the inverse matrix $(I+tD_x\bfw)^{-1}$
 and taking the inner products with  $l_1^*, l_2^*$, one
obtains two scalar equations for $\xi_{tt}$ and $\sigma_{tt}$.

Combined with (\ref{ett}), these yield the system
\bel{sxtt1}
\left\{ \bega{rl}
\eta_t&=~H,\\[3mm]
\sigma_{tt}&=~l_1^*\cdot  (I+t D_x\bfw )^{-1}\bigl[ F  + \eta^{-1} G\bigr],\\[3mm]
\xi_{tt}&\ds=~
{1\over 1+ t \lambda_2}\, l_2^*\cdot
\bigl[ F  + \eta^{-1}  G\bigr],
\enda\right.\eeq
where $H,F,G$ are the functions introduced at (\ref{ett})-(\ref{ytt}).
 Note the the second equation in (\ref{sxtt1}) can be equivalently written as
\bel{sxtt2}\sigma_{tt}~=~\frac{1}{1+t\lambda_1}\,l_1\cdot\bigl[ F  + \eta^{-1} G\bigr]+\frac{l_1^*\cdot r_2}{1+t\lambda_2}\,l_2\cdot\bigl[ F  + \eta^{-1} G\bigr].\eeq

\v

 {\bf 3.} Our next goal is to
 write the system (\ref{sxtt1}) in a more standard form, to which
a general existence-uniqueness theorem can then be applied. Here the main difficulty lies in the singularities along the initial curve where
$t=\tau(\zeta)$.   Indeed, the functions $F,G,H$ are not analytic, because some terms become unbounded as $t\to \tau(\zeta)$.  Moreover, as
shown in Lemma~\ref{l:42}, the right hand sides of (\ref{vpm2})-(\ref{rpm2})
 contain a square root  as a factor.  To motivate our future choices
 of rescaled variables, it is useful to keep in mind the order of magnitude of the various quantities:
$$\eta~=~\O(1)\cdot \sqrt{t-\tau(\zeta)}, \quad \qquad y_t - \bfv^\pm~=~\O(1)\cdot \sqrt{t-\tau(\zeta)},\qquad \quad \rho^\pm~=~{\O(1)\over t-\tau(\zeta)} \,.$$

We will write the equations (\ref{sxtt1}) in a form where all the singularities  arise
in connection with negative powers of $t-\tau(\zeta)$.
Toward this goal, we consider the rescaled variables
\bel{resv}
\Hat\eta\,=\, {\eta\over\sqrt{t-\tau(\zeta)}}\,,
\qquad\quad \Hat \xi\,=\, {\xi-\zeta\over t-\tau(\zeta)},
\qquad\quad  \Hat\sigma\,=\, {\sigma\over t-\tau(\zeta)}\,.\eeq
In addition, to reduce the second order equations in (\ref{sxtt1}) to a first order system,
we introduce the variables
\bel{SZdef} S\,=\, \sigma_t\,,\qquad\qquad Z\,=\, \xi_t\,.\eeq

\begin{remark} {\rm The results in Lemma~\ref{l:42} describe the singularities of the functions $\bfv^\pm, \rho^\pm$
in terms of the variables $(t,\xi,\Tilde\sigma)$.
However, we ultimately seek a solution in the form
$(\xi,\sigma,\eta)=(\xi,\sigma,\eta)(t,\zeta)$.  Toward this goal, it is convenient to express all the singular terms
appearing
in the equations as functions of the independent variables $t,\zeta$   (instead of $t,\xi$).

Notice the difference between $\Hat\sigma$ and the variable $\Tilde\sigma$ introduced at (\ref{Tsig}).
Namely,
$$\Tilde \sigma~=~\frac{t-\tau(\zeta)}{t-\tau\bigl(\xi(t,\zeta)\bigr)}\cdot\Hat \sigma\,.$$
When $t=\tau(\zeta)$, the boundary condition yields $\xi\bigl(\tau(\zeta), \zeta\bigr)=\zeta$.
Here the right hand side is analytic w.r.t.~the variables $(t, \zeta, \xi, \Hat \sigma)$, and hence also w.r.t.~$(t, \zeta,\Hat \xi, \Hat \sigma)$}
\end{remark}

\v
{\bf 4.}
We now analyze the structure of $H, F$ and $G$ in (\ref{sxtt1}).
To simplify the notation, it will be convenient to write
\bel{HB1}\Hat \bfB_1~\doteq~\bigl(t-\tau(\xi)\bigr)\bfB_1\,.\eeq
Recalling the notations and results in Lemma \ref{l:42}, we have
{\small
\begin{eqnarray*}H&=&\frac{1}{\sqrt{t-\tau(\xi)}}\cdot\bigg\{ -2 \bfB_2b_1-2 \Hat\bfB_1b_2+2 b_2
 (I+t D_x\bfw )\cdot  (x_\xi \xi_t + x_\sigma \sigma_t)\bigg\}\times \Big[(I+t D_x\bfw )\cdot  (x_\xi \xi_\zeta + x_\sigma \sigma_\zeta)\Big]\\[3mm]
 &\doteq&\frac{1}{\sqrt{t-\tau(\zeta)}}H_1(t,\zeta,\Hat \xi, \Hat\sigma, \xi_t,\sigma_t, \xi_{\zeta}, \sigma_{\zeta}), \end{eqnarray*}
 }
where
\bel{H1def} H_1=\sqrt{\frac{t-\tau(\zeta)}{t-\tau(\xi)}}\bigg\{ -2 \bfB_2b_1-2 \Hat\bfB_1b_2+2 b_2
 (I+t D_x\bfw )\cdot  (x_\xi \xi_t + x_\sigma \sigma_t)\bigg\}\times \Big[(I+t D_x\bfw )\cdot  (x_\xi \xi_\zeta + x_\sigma \sigma_\zeta)\Big].\eeq

 Observing that
$$\Hat \eta_t~=~{\eta_t\over \sqrt{t-\tau(\zeta)}}-{\Hat \eta\over 2\bigl(t-\tau(\zeta)\bigr)},$$
we obtain
\bel{H2def}\Hat\eta_t\,=\,\frac{H_2}{t-\tau(\zeta)},\qquad\qquad
 H_2\,\doteq \,H_1-\frac{\Hat\eta}{2}.\eeq

Next, we observe  that $F$ in (\ref{ytt}) is an analytic function of the variables
$(t,\zeta,\Hat \xi, \Hat\sigma, \xi_t,\sigma_t)$.
In particular, for a solution of the form
$$\sigma(t)~=~S_0 \bigl(t-\tau(\xi)\bigr) + \O(1)\cdot \bigl(t-\tau(\xi)\bigr)^2,
\qquad\qquad \xi(t)~=~\xi_0+ \O(1)\cdot \bigl(t-\tau(\xi)\bigr)^2,$$
we have
$$F~=~2\omega_1(\xi_0) S_0  r_1(\xi_0)+ \O(1)\cdot \bigl(t-\tau(\xi)\bigr).$$

By the first two equations in (\ref{yte}) and by Lemma \ref{l:42} it follows
\bel{Gcom}
\left\{
\bega{rl} y_t-\bfw&=~ (I+t  D_x \bfw)(x_\xi\xi_t + x_\sigma \sigma_t)~\doteq~{\bf A}\,,\\[3mm]
 \bfv^\pm-\bfw&=~\bigl(t-\tau(\xi)\bigr)\bfB_1 \pm \sqrt{t-\tau(\xi)}\,\bfB_2 , \\[3mm]
y_\zeta&=~(I+tD_x\bfw) (x_\xi \xi_\zeta+x_\sigma \sigma_\zeta),\\[4mm]
\rho^\pm&=\ds~{b_1\pm b_2\sqrt{t-\tau(\xi)}\over t-\tau(\xi)}\,.\enda\right.\eeq

 To leading order, i.e., neglecting terms of order $\O(1)\cdot \bigl(t-\tau(\xi)\bigr)^{3/2}$,
 recalling (\ref{bet0}) and (\ref{bb}) we find
\bel{app3}\bfA ~\approx~\bigl(1+ t\lambda_1(\xi)\bigr) r_1(\xi) S_0~=~-
\bigl(t-\tau(\xi)\bigr) \omega_1(\xi) r_1(\xi) S_0\,,\eeq
\bel{app4}\bega{rl} \bfv^\pm-\bfw&\approx~\lambda_1(\xi) \bigl[ Z^\pm-\sigma\bigr]r_1(\xi)\\[2mm]
&\approx~-\omega_1(\xi)\left[   \bigl(t-\tau(\xi)\bigr)\cdot\beta_0(t,\xi,\Tilde\sigma)\, \pm\, \sqrt{ t-\tau(\xi)} \cdot \beta_1(t,\xi,\Tilde\sigma) - S_0 \bigl(t-\tau(\xi)\bigr)\right] \\[2mm]
&\ds\approx~-\omega_1(\xi)\left[  - {3\over 2} S_0  \bigl(t-\tau(\xi)\bigr) \pm \omega_1\sqrt{ 6\bigl(t-\tau(\xi)\bigr)\over
\omega_3}- \Tilde\phi(\xi) \cdot \bigl(t-\tau(\xi)\bigr) \right] r_1(\xi).
\enda \eeq
Here $\Tilde\phi$ collects terms that do not depend on $\sigma$.
Notice that here $\bfw=\bfw\bigl(x(\xi(t), \sigma(t))\bigr)$ must be computed along a solution of (\ref{sxtt1}).
Similarly, the values $Z^\pm$ also refer to this solution.  For this reason, they are given by
(\ref{Zrep}), with $\Tilde \sigma= S_0$.

Combining (\ref{app3})-(\ref{app4}), we thus reach the expression
\bel{app5}\bega{rl} \bfv^\pm - y_t&=~\ds
\omega_1(\xi)\left[ {5\over 2} S_0 \bigl(t-\tau(\xi)\bigr)\mp \omega_1 (\xi)\sqrt{6(t-\tau(\xi))\over \omega_3(\xi)}
+ \Tilde\phi(\xi) \cdot \bigl(t-\tau(\xi)\bigr)\right] r_1(\xi)\\[4mm]
&\qquad+\O(1)\cdot \bigl(t-\tau(\xi)\bigr)^{3/2}.\enda\eeq

We note that
\bel{sizes}
\Hat\bfB_1\,=\, \O(1)\cdot \bigl(t-\tau(\xi)\bigr),\qquad
 \bfB_2~=~\O(1),\qquad {\bf A}\times y_\zeta\,=\, \O(1)\cdot \bigl(t-\tau(\xi)\bigr),\eeq
%
 It then follows that
$$l_1^*\cdot(I+t D_x\bfw )^{-1}G~=~ \O(1)\cdot \bigl(t-\tau(\xi)\bigr)^{-1/2}.$$
In addition we have
$$l_1^*\cdot(I+t D_x\bfw )^{-1}F~=~\O(1)\cdot \bigl(t-\tau(\xi)\bigr)^{-1}.$$
%

Since $l_1^*\cdot(I+tD_x{\bfw})^{-1}=\O(1)\cdot(t-\tau(\xi))^{-1}$,  recalling (\ref{sxtt1})
and (\ref{SZdef}), we have
\bel{hxi}\sigma_{tt}~=~S_t\,=\,\frac{F_1}{t-\tau(\zeta)},\eeq
where
\bel{F1def} F_1\,\doteq\,(t-\tau(\zeta))l_1^*\cdot(I+tD_x{\bfw})^{-1}F+\sqrt{\frac{t-\tau(\zeta)}{t-\tau(\xi)}}\cdot\frac{1}{\widehat{\eta}}\,l_1^*\cdot(I+tD_x{\bfw})^{-1}G_0\,.\eeq
 \bel{G0def}G_0\,\doteq\, \sqrt{t-\tau(\xi)}\cdot G.\eeq

In a similar way, the third equation in (\ref{sxtt1}) now becomes
\bel{hsi}Z_t~=~\left(l_2^*\cdot r_2\right)\frac{l_2\cdot F}{1+t\lambda_2}+\left(l_2^*\cdot r_2\right)\frac{l_2\cdot G}{1+t\lambda_2}\frac{1}{\widehat{\eta}\sqrt{t-\tau(\zeta)}}~\doteq ~F_2+\frac{G_2}{t-\tau(\zeta)},\eeq
where \bel{F2G2}
F_2~=~\left(l_2^*\cdot r_2\right)\frac{l_2\cdot F}{1+t\lambda_2},\qquad
\qquad G_2~=~\sqrt{\frac{t-\tau(\zeta)}{t-\tau(\xi)}}\,\frac{l_2^*\cdot r_2}{1+t\lambda_2}\,\frac{1}{\widehat{\eta}}\, l_2\cdot G_0\,,\eeq are analytic functions
of their arguments.

\v
{\bf 5.}
 Summarizing the above discussion, the second order system (\ref{sxtt1}) can be rewritten
 as a first order system of five  equations. Namely
\bel{etat}
\Hat\eta_t~=~\frac{H_2}{t-\tau(\zeta)},\qquad \qquad H_2\doteq H_1-\frac{\Hat\eta}{2} ,\eeq
for a suitable analytic function $ H_1$, together with
\bel{sit}
\Hat\sigma_t~=~{1\over t-\tau(\zeta)} (S-\Hat \sigma),\eeq
\bel{XIT}\Hat \xi_t~=~{1\over t-\tau(\zeta)} (Z-\Hat \xi),\eeq
\bel{Hxi}
S_t~=~{1\over t-\tau(\zeta)} F_1 \bigl( t,\zeta, \Hat \xi, \Hat \sigma,\Hat \eta,  \xi_t,\sigma_t, \xi_\zeta, \sigma_\zeta\bigr),\eeq
\bel{Hsi}
Z_t~=~F_2 \bigl( t,\zeta, \Hat \xi, \Hat \sigma, \Hat \eta,  \xi_t,\sigma_t, \xi_\zeta, \sigma_\zeta\bigr) + {1\over  t-\tau(\zeta)}
 G_2 \bigl( t,\zeta, \Hat \xi, \Hat \sigma, \xi_t,\sigma_t, \xi_\zeta, \sigma_\zeta\bigr) ,\eeq
for suitable analytic functions $F_1, F_2, G_2$.
Boundary conditions are assigned for $t=\tau(\zeta)$.
In addition to
\bel{IC0} \eta \bigl(\tau(\zeta), \zeta\bigr)~=~ \sigma\bigl(\tau(\zeta), \zeta\bigr)~=~0,\qquad\qquad  \xi\bigl(\tau(\zeta), \zeta\bigr)~=~\zeta,\eeq
we impose the identities
\bel{IC1} \left\{ \bega{rll} \Hat\eta\bigl(\tau(\zeta), \zeta\bigr)&=~\varrho_0(\zeta),\\[3mm]
\Hat\sigma\bigl(\tau(\zeta), \zeta\bigr)&=~S\bigl(\tau(\zeta), \zeta\bigr)&=~S_0(\zeta),\\[3mm]
\Hat\xi\bigl(\tau(\zeta), \zeta\bigr)&=~Z\bigl(\tau(\zeta), \zeta\bigr)&=~Z_0(\zeta),
\enda\right.\eeq
where the three functions  $\varrho_0, S_0, Z_0$ are uniquely determined by the requirement
that $H_2, F_1, G_2$ in (\ref{etat}), (\ref{Hxi}), (\ref{Hsi}) vanish at $t=\tau(\zeta)$.
\v
{\bf 6.} In this step
we show how to uniquely determine the initial data (\ref{IC1}).

We begin by showing that $Z_0(\zeta)=0$.
At $t=\tau(\zeta)$, by (\ref{F2G2})  the requirement $G_2=0$ is equivalent to $l_2\cdot G_0=0$.
Recalling the bounds
(\ref{sizes}), from  $l_2\cdot G_0=0$ we obtain
\bel{Z01}\big( \bfB_2\times y_\zeta\big)\big(l_2\cdot(\Hat\bfB_1-{\bf A})\big)=0.\eeq
We claim that $\bfB_2\times y_\zeta \neq0$. Indeed,  at $t=\tau(\zeta)$ one has $\xi_{\zeta}=1$, and hence  by our assumptions it follows
\bel{Z02}\bega{rl} \bfB_2\times y_\zeta&\ds=~\bfB_2\times \Big[(I+\tau(\zeta)D_x\bfw)\cdot r_2^*(\zeta)\Big]~=~\frac{1}{\lambda_1(\zeta)}\bfB_2\times
\Big[\bigl(\lambda_1(\zeta)I-D_x\bfw\bigr)\cdot r_2^*(\zeta)\Big]\\[3mm]
&=\ds~\sqrt{c_1(\zeta)}\, r_1(\zeta)\times \Big[\bigl(\lambda_1(\zeta)I-D_x\bfw\bigr)\cdot r_2^*(\zeta)\Big],\enda\eeq
which does not vanish, since $r_2^*$ is not parallel to $r_1$.
Here $c_1(\cdot)$ is the function introduced at (\ref{bb}).
Moreover,  at $t=\tau(\zeta)$ we have
\bel{Z03}\bega{rl}l_2\cdot (\Hat\bfB_1-{\bf A})&=\,\ds-l_2\cdot(I+tD_x\bfw)r_2^*(\zeta)Z_0(\zeta)\\[3mm]
&=\,\ds-l_2\cdot\bigl(I-\frac{1}{\lambda_1(\zeta)}D_x\bfw\bigr)r_2^*(\zeta)Z_0(\zeta)\\[3mm]
&=\,\ds-\frac{1}{\lambda_1(\zeta)}l_2\cdot\bigl(\lambda_1(\zeta)I-D_x\bfw\bigr)r_2^*(\zeta)Z_0(\zeta)\\[3mm]
&=\,\ds-\frac{1}{\lambda_1(\zeta)}\bigl(\lambda_1(\zeta)-\lambda_2(\zeta)\bigr)l_2\cdot r_2^*(\zeta)Z_0(\zeta),\enda\eeq
because $l_2$ is a left eigenvector of $D_x\bfw$.
Combining (\ref{Z01})--(\ref{Z03}), we conclude that $Z_0(\zeta)=0$.
\v

 Next, we show how to determine $\varrho_0(\zeta)$.  All computations are performed in the limit $t\to \tau(\zeta)$. Notice that in this case $\xi=\zeta$, hence  $\tau(\xi)=\tau(\zeta)$.

 The requirement $H_2=0$ implies
$\varrho_0(\zeta)=2H_1.$
By  (\ref{HB1}) it trivially follows $\Hat\bfB_1=\O(1)\cdot(t-\tau(\xi))$.
Since $Z_0(\zeta)=0$, by (\ref{bB0}) at $t=\tau(\zeta)$ we obtain
\bel{rho>0}\bega{rl} \varrho_0(\zeta)&=~2H_1~=\,-4b_1\,\bfB_2\times \Big[(I+\tau(\zeta)D_x\bfw)r_2^*(\zeta)\Big]\\[4mm]
&~=\,-4b_1\sqrt{c_1(\zeta)}\, r_1(\zeta)\times \Big[(\lambda_1I-D_x\bfw)r_2^*(\zeta)\Big] ~>~0.\enda\eeq
Here the inequality $\varrho_0>0$ is intuitively obvious, because the density
along the singular curve is initially zero and positive afterwards.   A more detailed proof
is achieved by observing that, at $t=\tau(\zeta)$, one has
$$r_2^*(\zeta)~=~a_1 r_1(\zeta) + a_2 r_2(\zeta),$$ with $a_2>0$.
Recalling that $b_1$
is the positive constant at (\ref{b1e}), the right hand side of (\ref{rho>0}) can be written as
$$-4b_1\sqrt{c_1(\zeta)}\, r_1(\zeta)\times \Big[(\lambda_1I-D_x\bfw) a_2 r_2(\zeta)\Big]
~=\,-4b_1\sqrt{c_1(\zeta)}\, r_1(\zeta)\times \Big[(\lambda_1-\lambda_2) a_2 r_2(\zeta)\Big] ~>~0.$$

Finally, in order to determine the initial value $S_0(\zeta)$, we need to use the assumption $F_1=0$ at $t=\tau(\zeta)$.
Since $\widehat{\eta}=\varrho_0(\zeta)$ at $t=\tau(\zeta)$,  one only needs to consider the equation
\bel{F1eq}F_1~=~(t-\tau(\zeta))l_1^*\cdot(I+tD_x{\bfw})^{-1}F+\frac{1}{\varrho_0(\zeta)}\,l_1^*\cdot(I+tD_x{\bfw})^{-1}G_0~=~0\,.\eeq
Inserting the expression (\ref{ytt}) for $F$ into the above equation and recalling that $Z_0(\zeta)=0$, always for $t=\tau(\zeta)$ we obtain
\bel{p1} \bega{l}\ds (t-\tau(\zeta))l_1^*\cdot(I+tD_x{\bfw})^{-1}F\\[4mm]
\quad \ds =(t-\tau(\zeta))l_1^*\cdot(I+tD_x{\bfw})^{-1}\Big\{- 2D_x\bfw \cdot  (x_\xi \xi_t + x_\sigma \sigma_t) \\[4mm]
\qquad\qquad  -
(t\, D^2\bfw ) \bigl((x_\xi \xi_t + x_\sigma \sigma_t)
 \otimes  (x_\xi \xi_t + x_\sigma \sigma_t) \bigr)
  -(I+tD_x\bfw)(x_{\xi\xi}\xi_t^2+2x_{\xi\sigma}\xi_t\sigma_t+x_{\sigma\sigma}\sigma_t^2)\Big\}\\[4mm]
\quad \ds =(t-\tau(\zeta))l_1^*\cdot(I+tD_x{\bfw})^{-1}\Big\{- 2D_x\bfw \cdot r_1(\zeta)S_0(\zeta) -\tau(\zeta)
 S_0^2(\zeta)D^2\bfw(\zeta)\bigl(r_1(\zeta)\otimes r_1(\zeta)\bigr)\Big\} \\[4mm]
\quad \ds =(t-\tau(\zeta))l_1^*\cdot(I+tD_x{\bfw})^{-1}\Big\{- 2\lambda_1(\zeta)r_1(\zeta)S_0(\zeta) -\tau(\zeta)
 S_0^2(\zeta)D^2\bfw(\zeta)\bigl(r_1(\zeta)\otimes r_1(\zeta)\bigr) \Big\}\\[4mm]
\quad \ds =~  -2S_0(\zeta)-S_0^2(\zeta)\Big[\tau(\zeta)(t-\tau(\zeta))l_1^*\cdot(I+tD_x{\bfw})^{-1}
 r_1^T(\zeta)\, D^2\bfw \,r_1(\zeta)\Big],\enda\eeq
 where we have used the identity
\bel{di1} (t-\tau(\zeta))l_1^*\cdot(I+tD_x{\bfw})^{-1}r_1~=~\frac{1}{\lambda_1(\zeta)}.\eeq
At first sight, (\ref{p1}) looks like a quadratic equation for $S_0$.
However, the following analysis reveals that as $t\to\tau(\zeta)$
the coefficient of $S_0^2$ vanishes, while the
coefficient of $S_0$ is nonzero.  Therefore, a unique value of $S_0(\zeta)$ is determined.

Indeed, consider the quantity
\bel{te11}
l_1^*\cdot (I+tD_x\bfw)^{-1} (t D^2\bfw)(r_1\otimes r_1).\eeq
For any vector $\bfv$, decomposing $\bfv$ along the basis of eigenvectors
$\{r_1,r_2\}$, one finds
$$(I+tD_x)\bfv~=~(l_1\cdot \bfv) (1+t\lambda_1) r_1 + (l_2\cdot \bfv) (1+t\lambda_2) r_2\,,$$
$$(I+tD_x)^{-1} \bfv~=~(l_1\cdot \bfv) {1\over 1+t\lambda_1} r_1 + (l_2\cdot \bfv) {1\over 1+t\lambda_2} r_2\,.$$
Therefore, the quantity in (\ref{te11}) can be rewritten as
\bel{te22} l_1^* \cdot \left[
 {1\over 1+t\lambda_1}  \Big( l_1\cdot (t D^2\bfw)(r_1\otimes r_1)\Big) r_1+
  {1\over 1+t\lambda_2}  \Big( l_2\cdot (t D^2\bfw)(r_1\otimes r_1)\Big) r_2\right].\eeq
For $t\approx -1/\lambda_1$ the second term remains bounded, but the first term yields
$$
 {1\over 1+t\lambda_1}  \Big( l_1\cdot (t D^2\bfw)(r_1\otimes r_1)\Big)( l_1^* \cdot  r_1)
 ~=~ {1\over 1+t\lambda_1}  \Big( l_1\cdot (t D^2\bfw)(r_1\otimes r_1)\Big).
$$
The following lemma allows us  to compute the quantity
$ l_1\cdot (D^2\bfw)(r_1\otimes r_1)$.

\begin{lemma} \label{l:sder}
Let $\bfw(x)$ be a smooth vector field on $\R^2$.  Assume that  at each point
$x\in\R^2$ its Jacobian matrix
$D_x\bfw(x)$ has eigenvalues $\lambda_1(x)< \lambda_2(x)$, and a basis of right eigenvectors $\bigl\{ r_1(x), r_2(x)\bigr\}$. Call $\bigl\{ l_1(x), l_2(x)\bigr\}$
the dual basis of left eigenvectors.

Assume that at a point $\bar x$ we have $\nabla\lambda_1(\bar x)\cdot r_1(\bar x)=0$.
Then
\bel{sderw} l_1(\bar x) \cdot\Big[ \bigl(D^2\bfw(\bar x)\bigr)\bigl(r_1(\bar x)\otimes r_1(\bar x)\bigr)\Big]~=~0.\eeq
\end{lemma}

{\bf Proof.} For any vector $\bfv\in\R^2$, the directional derivative of $\bfw$
in the direction of $\bfv$ at a point $x\in\R^2$ is computed by
$${d\over ds} \bfw(x+s\bfv)\Bigg|_{s=0}~=~
\lambda_1(x) \bigl(l_1(x)\cdot \bfv\bigr) r_1(x) + \lambda_2(x) \bigl(l_2(x)\cdot \bfv\bigr) r_2(x).$$
Calling $\bar r_1\doteq r_1(\bar x)$, we need to compute the second derivative:
\bel{q1}\bega{l}\ds
\bigl(D^2\bfw(\bar x)\bigr)\bigl(r_1(\bar x)\otimes r_1(\bar x)\bigr)~=~ {d^2\over ds^2} \bfw(\bar x+s\bar r_1)\Bigg|_{s=0}\\[4mm]
 \ds =~{d\over ds}\Big[
\lambda_1(\bar x+ s \bar r_1 ) \bigl(l_1(\bar x+ s \bar r_1 )\cdot \bar r_1\bigr) r_1(\bar x+ s \bar r_1 ) + \lambda_2(\bar x+ s \bar r_1 ) \bigl(l_2(\bar x+ s \bar r_1 )\cdot \bar r_1\bigr) r_2(\bar x+ s \bar r_1 )\Big]_{s=0}
\enda\eeq
Since
$$l_1(\bar x) \cdot r_1(\bar x)\,=\,1,\qquad\qquad l_1(\bar x) \cdot r_2(\bar x)\,=\,0,
\qquad\qquad
\nabla\lambda(\bar x)\cdot r_1(\bar x)=0,$$
taking the product of the right hand side of (\ref{q1}) with $l_1(\bar x)$ some terms
clearly vanish.  The remaining terms are computed by
\bel{q2} \bega{l}
 l_1(\bar x) \cdot\Big[ \bigl(D^2\bfw(\bar x)\bigr)\bigl(r_1(\bar x)\otimes r_1(\bar x)\bigr)\Big]\\[3mm]
\ds \quad =~\lambda_1(\bar x) \left[ {d\over ds} \bigl(l_1(\bar x+ s \bar r_1 )\cdot \bar r_1\bigr) + l_1(\bar x)  \cdot  {d\over ds} r_1(\bar x+ s\bar r_1)\right]_{s=0}\\[4mm]
\ds  \quad =~\lambda_1(\bar x) \Big[ {d\over ds} \bigl(l_1(\bar x+ s \bar r_1 )\cdot  r_1(\bar x+ s \bar r_1 )\bigr)\Big]_{s=0}~=~0.
 \enda\eeq
\endproof
By Lemma \ref{l:sder} together with (\ref{te11})-(\ref{te22}) we conclude that, on the right hand side of (\ref{p1}), the term containing $S_0^2$  vanishes  at $t=\tau(\zeta)$.
Notice that this is entirely analogous to the 1-D case, where the terms containing $S^2$ on the right hand side of (\ref{SO1}) vanish as $t\to t_0$.

By the previous analysis, the equation (\ref{F1eq}) further reduces to
\bel{F11}-2 S_0(\zeta) +\frac{1}{\varrho_0(\zeta)}\,l_1^*\cdot(I+tD_x{\bfw})^{-1}G_0~=~0\,.\eeq
We will show that (\ref{F11}) yields a linear equation for $S_0(\zeta)$,  where the coefficient of $S_0$ does not vanish. This will uniquely determine $S_0$.

Toward this goal,
note that $G_0=\sqrt{t-\tau(\xi)}\,G$, with
\bel{G22}\bega{rl} G&=~
\bigl( (\bfv^--y_t) \times y_\zeta\bigr) \rho^- (\bfv^- - y_t)-
\bigl( (\bfv^+-y_t)\times y_\zeta \bigr) \rho^+ (\bfv^+ - y_t).
\enda
\eeq
Therefore, in the expression of $G$ we should analyze the terms of order $\sqrt{t-\tau(\xi)}$ which also contain $S_0$, as $t\rightarrow \tau(\zeta)$.
Based on (\ref{app5}) and Lemma \ref{l:41}, to leading order we compute
 \bel{G23}\bega{rl} G&\ds\approx~\omega_1^2\left[\frac52 S_0(t-\tau(\xi))+\omega_1\sqrt{\frac{6(t-\tau(\xi))}{\omega_3}}+\Tilde\phi(\xi)(t-\tau(\xi))\right]^2\cdot \bigl(r_1\times y_{\zeta}\bigr)\rho^-\,r_1\\[3mm]
 &\ds\qquad -\omega_1^2\left[\frac52 S_0(t-\tau(\xi))-\omega_1\sqrt{\frac{6(t-\tau(\xi))}{\omega_3}}+\Tilde\phi(\xi)(t-\tau(\xi))\right]^2
 \cdot \bigl(r_1\times y_{\zeta}\bigr)\rho^+\,r_1\\[3mm]
 &\ds=~\sqrt{t-\tau(\xi)}\omega_1^2\bigl(r_1\times y_{\zeta}\bigr)r_1\Bigl\{10b_1S_0\omega_1\sqrt{\frac{6}{\omega_3}}-b_2\frac{12\omega_1^2}{\omega_3}\Bigr\}+\O(1)(t-\tau(\xi))^{3/2}.
\enda
\eeq
Hence,  using (\ref{G23}) and  (\ref{di1}) we obtain
\bel{G24} \bega{rl}\ds \frac{1}{\varrho_0(\zeta)}\,l_1^*\cdot(I+tD_x{\bfw})^{-1}G_0
&\ds=~\frac{\bigl(r_1\times y_{\zeta}\bigr)\omega_1^2}{-4b_1\sqrt{c_1(\zeta)}\;\lambda_1^2\bigl(r_1\times y_{\zeta}\bigr)}
\left[10b_1S_0\omega_1\sqrt{\frac{6}{\omega_3}}-b_2\frac{12\omega_1^2}{\omega_3}\right]
\\[4mm]
&=\ds -\frac{10}{4}S_0+\frac{b_2}{2b_1}\sqrt{c_1},\enda\eeq
where we have made use of the formula for $\varrho_0(\zeta)$  at (\ref{rho>0}).

We claim  that, at $t=\tau(\zeta)$,
 \bel{b21}\frac{b_2}{2b_1}\sqrt{c_1(\zeta)}~=~\frac34 S_0+\text{other terms independent of $S_0$}.\eeq

 Indeed, since
$$\lambda_1(X^{\pm})-\lambda_1(\xi)=\frac{\omega_3}{2}\left(Z^{\pm}\right)^2+\frac{\omega_4}{6}\left(Z^{\pm}\right)^3+\O(1)\left(Z^{\pm}\right)^4,$$
by Lemma \ref{l:41} and (\ref{be0})-(\ref{be1}) we have
$$\lambda_1(Z^{\pm})-\lambda_1(\xi)~=~\bigl(t-\tau(\xi)\bigr)\left\{\frac{\omega_3}{2}c_1(\xi)+\O(1)\bigl(t-\tau(\xi)\bigr)\right\}
\pm\sqrt{t-\tau(\xi)}\Bigl[P+\O(1)\bigl(t-\tau(\xi)\bigr)\Bigr],$$
where $$P~\doteq~\omega_3\sqrt{c_1(\xi)}\left(-{c_2(\xi,\Tilde\sigma)+c_1(\xi)c_3(\xi)\over 2}\right)+\frac{\omega_4}{6}c_1(\xi)^{3/2}.$$
Therefore,
\begin{eqnarray*}1+t\lambda_1(\xi,\sigma)&=&\frac{2}{\tau(\xi)}\bigl(t-\tau(\xi)\bigr)\left[1+\O(1)\bigl(t-\tau(\xi)\bigr)
\pm\sqrt{t-\tau(\xi)}\Bigl(\frac{\tau(\xi)^2}{2}P+\O(1)\bigl(t-\tau(\xi)\bigr)\Bigr)\right]\\[4mm]
&\doteq&\frac{2}{\tau(\xi)}\bigl(t-\tau(\xi)\bigr)\left[1+\phi_3\pm\sqrt{t-\tau(\xi)}\;\phi_4\right],\end{eqnarray*}
where $\phi_3, \phi_4$ are  analytic functions.
It further implies
\bel{mn}\frac{1}{1+t\lambda_1}=\frac{\tau(\xi)}{2\bigl(t-\tau(\xi)\bigr)}\Bigl[1\pm\sqrt{t-\tau(\xi)}\;\vp_0\Bigr],\eeq
where $\vp_0$ is an analytic function satisfying
$$t=\tau(\xi),\qquad \vp_0=-\frac{\tau(\xi)^2}{2}\left(\omega_3\sqrt{c_1(\xi)}\,\Lambda+\frac{\omega_4}{6}c_1(\xi)^{3/2}\right).$$

Next, by (\ref{te1}) and (\ref{mn}) one obtains
\bel{b1b2}b_1=\frac{\tau(\xi)}{2}\left[\varphi_1 +(t-\tau(\xi))\varphi_2\vp_0\right],\qquad
b_2=\frac{\tau(\xi)}{2}\left[\varphi_1 \vp_0+\varphi_2\right].\eeq
Finally,  at $t=\tau(\zeta)$ we obtain
\begin{eqnarray*}\frac{b_2}{2b_1}\sqrt{c_1(\xi)}&=&\frac{(\varphi_1 \vp_0+\varphi_2)\sqrt{c_1(\xi)}}{2[\varphi_1 +(t-\tau(\xi))\varphi_2\vp_0]}\\[4mm]
&=&\sqrt{c_1(\zeta)}\left(\frac{\vp_0}{2}+\frac{\varphi_2}{\varphi_1}\right)~=~-\frac{\tau(\zeta)^2}{4}c_1(\zeta)\omega_3\,\Lambda+\sqrt{c_1(\zeta)}\frac{\varphi_2}{\varphi_1}\\[4mm]
&=&\frac{3}{4}S_0+\sqrt{c_1(\zeta)}\frac{\varphi_2}{\varphi_1},
\end{eqnarray*}
as claimed in (\ref{b21}). We note that $\sqrt{c_1(\zeta)}\frac{\varphi_2}{\varphi_1}$ does not depend on $\Tilde\sigma$ at $t=\tau(\zeta)$.

Notice that the above result is consistent with what found in the 1-D case.

 Thus, by combining (\ref{b21}) and (\ref{G24}) we conclude that the coefficient of $S_0$ in (\ref{F11}) is exactly $\frac{15}{4}$ as in 1-D case. Therefore,  from the equation $F_1=0$ at $t=\tau(\zeta)$
a unique value of $S_0(\zeta)$ is determined.

\v
 {\bf 7.} We now perform a further renaming of variables, using
 $$T\,=\,t-\tau(\zeta)\qquad\hbox{and}\qquad \zeta$$
 as independent variables, and setting
 $\bfu\doteq\bigl( \Hat \eta, \Hat\sigma, \Hat\xi, S,Z)\in \R^5$.
In this way, the original set of equations
 (\ref{etat})--(\ref{Hsi}) can be rewritten as a system of the form
\bel{SS}
\bfu_T~=~ {1\over T}\Phi(T, \zeta,\bfu,\bfu_{\zeta}) + \Psi(T, \zeta,\bfu, \bfu_{\zeta}).\eeq
The initial data, given  at $T=0$, are
\bel{SIN}\bfu(0, \zeta)~=~\ov \bfu(\zeta)~=~ \Big(\rho_0(\zeta),\,0,\,\zeta,\,S_0(\zeta),\, 0\Big).\eeq
where $\rho_0$ and $S_0$ are the functions considered at
(\ref{rho>0}) and (\ref{F1eq}), respectively.

A unique local solution of the Cauchy problem is sought, relying on the fact
that $\Phi, \Psi$ are analytic functions of their arguments, and moreover
\bel{ICO}\Phi\bigl(0, \zeta,\ov\bfu(\zeta),\ov\bfu_{\zeta}(\zeta)\bigr)~=~0\qquad\qquad\forall \zeta.\eeq


 After  changing the variables $\widehat{\bfu}=\bfu-\ov{\bfu}(\zeta)$ and  removing
 the  hats, we thus obtain the following Cauchy problem
\bel{SS3} \bfu_T~=~{1\over T} \bfF(T, \zeta, \bfu, \bfu_{\zeta}) + \bfG(T,\zeta,\bfu,\bfu_{\zeta}) ,\eeq
with initial data \bel{INC} \bfu(0,\zeta)=0,\eeq
where $\bfF,\bfG$ are analytic functions of their arguments and moreover
$$\bfF\bigl(0,\zeta,0,0\bigr)~=~0\qquad\qquad\forall \zeta.$$

\v
 {\bf 8.} Note that (\ref{SS3}) can be reformulated as a system of singular nonlinear PDEs
\bel{sip}T\,\bfu_T~=~\bfH(T,\zeta,\bfu, \bfu_{\zeta}),\eeq
where $\bfH(T,\zeta,\bfu,\bfv)=\bfF(T,\zeta,\bfu,\bfv)+T\cdot\bfG(T,\zeta,\bfu,\bfv)$, which is analytic of its arguments. The system (\ref{sip}) is called be of Fuchs type \cite{bao} or
Briot-Bouquet type  with respect to $T$ if
$$\bfH(0,\zeta,0,0)\,=\,0,\qquad \qquad \bfH_{\bfv}(0,\zeta,0,0)\,=\,0,$$
in a neighborhood of $\zeta=0$. The first condition is certainly satisfied. We need to check the condition $\bfF_{\bfv}(0,\zeta,0,0)=0$, which is equivalent to
\bel{ze}\frac{\partial \Phi}{\partial \bfv}\bigl(0,\zeta,\ov\bfu(\zeta),\ov\bfu'(\zeta)\bigr)~=~0.\eeq

 By (\ref{etat})-(\ref{Hsi}), we note that
$$\Phi~=~\Big(H_2, ~S-\Hat\sigma, Z-\Hat\xi, ~F_1, G_2\Big).$$
Then the condition (\ref{ze}) requires that, when $T=0$,
\bel{zefg}\frac{\partial H_2}{\partial \Hat\xi_{\zeta}}~=~ \frac{\partial H_2}{\partial \Hat\sigma_{\zeta}}~=~\frac{\partial F_1}{\partial \Hat\xi_{\zeta}}~=~\frac{\partial F_1}{\partial \Hat\sigma_{\zeta}}~=~ \frac{\partial G_2}{\partial \Hat\xi_{\zeta}}~=~ \frac{\partial G_2}{\partial \Hat\sigma_{\zeta}}~=~0.\eeq
Since at $T=0$ one has $\frac{\partial\xi_{\zeta}}{\partial\Hat\xi_{\zeta}}=0$, $
\frac{\partial\sigma_{\zeta}}{\partial\Hat\sigma_{\zeta}}=0$, we conclude that all the quantities in (\ref{zefg}) vanish, which in turn implies that (\ref{ze}) is valid. In other words, the system (\ref{sip}) is of Fuchs type or Briot-Bouquet type.

\v
 {\bf 9.}
Next, consider the Jacobian matrix
\bel{lad}{\bf\Lambda}(\zeta)~\doteq~\frac{\partial \bfH}{\partial \bfu}(0,\zeta,0,0).\eeq
By the previous analysis, we obtain
\bel{chm}{\bf\Lambda}(\zeta)~=~\frac{\partial\Phi}{\partial\bfu}\Big(0,\zeta,\ov\bfu(\zeta),\ov\bfu'(\zeta)\Big)~=~\left(\begin{array}{ccccc}-\frac12
&0&\frac{\partial H_2}{\partial\Hat\xi}&0&\frac{\partial H_2}{\partial Z}\\[2mm]
0&-1&0&1&0\\[2mm]
0&0&-1&0&1\\[2mm]
0&0&\frac{\partial F_1}{\partial\Hat\xi}
&\frac{\partial F_1}{\partial S}&\frac{\partial F_1}{\partial Z}\\[2mm]
0&0&0
&0&0\end{array}\right).\eeq
The eigenvalues of ${\bf\Lambda}(\zeta)$ are found by solving
\bel{eig}\lambda\left(\lambda+\frac12\right)(\lambda+1)^2\left(\lambda-\frac{\partial F_1}{\partial S}\right)\,=\,0,\eeq
where $\frac{\partial F_1}{\partial S}$ is evaluated at $t=\tau(\zeta)$.
We have
$$F_1~=~(t-\tau(\zeta))l_1^*\cdot(I+tD_x{\bfw})^{-1}F+\frac{1}{\widehat{\eta}}\sqrt{\frac{t-\tau(\zeta)}{t-\tau(\xi)}}{}\,l_1^*\cdot(I+tD_x{\bfw})^{-1}G_0\,,$$
where  $$G_0~=~-2(t-\tau(\xi))b_2({\bfB_2}\times y_\zeta ){\bfB_2}-2b_1({\bfB_2}\times y_\zeta)N+2b_2(M\times y_\zeta)N+2b_1(M\times y_\zeta){\bfB_2}.$$ We note that the variable $S$ only appears in $F$, $M$ and $N$. Thus, at $T=0$ we have
$$\frac{\partial F_1}{\partial S}=(t-\tau(\zeta))l_1^*\cdot(I+tD_x{\bfw})^{-1}\frac{\partial F}{\partial S}
+\frac{1}{\varrho_0(\zeta)}\,l_1^*\cdot(I+tD_x{\bfw})^{-1}\frac{\partial G_0}{\partial S}.$$
For convenience, we introduce the notation
$${\bf d}(t,\xi,\sigma)~=~(t-\tau(\zeta))l_1^*\cdot(I+tD_x{\bfw})^{-1},$$
so that ${\bf d}(\tau(\zeta),\zeta,0)\doteq {\bf d}(\zeta)\neq0$. Indeed, together with (\ref{di1})
we have the identity
\bel{di}l_1^*\cdot(I+tD_x{\bfw})^{-1}r_1=\frac{1}{1+t\lambda_1}\,.\eeq

At $T=t-\tau(\zeta)=0$, by (\ref{di1}) we compute
\bel{f1s}\bega{rl}\ds\frac{\partial F_1}{\partial S}&\ds=\,-2{\bf d}(\zeta)\left[\lambda_1(\zeta)r_1(\zeta)+\tau(\zeta)D^2{\bfw}\,\bigl(r_1(\zeta)\otimes r_1(\zeta)\bigr)S_0(\zeta)\right]\\[2mm]
&\qquad\qquad \ds
+\frac{1}{\varrho_0(\zeta)}\Big[4\sqrt{c_1(\zeta)}\,b_1\lambda_1(\zeta)\bigl(r_1\times y_{\zeta}\bigr)\Big]\\[4mm]
&\ds=\,-2+\frac{1}{\varrho_0(\zeta)}\Big[4\sqrt{c_1(\zeta)}\,b_1\lambda_1(\zeta)\bigl(r_1\times y_{\zeta}\bigr)\Big]~=~-3,\enda\eeq
where we have used the identity $\varrho_0(\zeta)=-4b_1\sqrt{c_1(\zeta)}\lambda_1(\zeta)\bigl(r_1\times y_\zeta\bigr)$ together with
 $${\bf d}(\zeta)\tau(\zeta)D^2{\bfw}\,r_1(\zeta)\otimes r_1(\zeta)~=~0,$$
as proved in Lemma \ref{l:sder}.

Here we work out the computations for
 $$l_1^*\cdot(I+tD_x{\bfw})^{-1}\frac{\partial G_0}{\partial S}~=~\frac{\partial}{\partial S}\left[l_1^*\cdot(I+tD_x{\bfw})^{-1} G_0\right].$$
We first consider
  \bel{lg0}\bega{l}\ds l_1^*\cdot(I+tD_x{\bfw})^{-1}G_0\\[3mm]
\qquad =~\ds 2 l_1^*\cdot(I+tD_x{\bfw})^{-1} \bigg\{  \bigl((\bfA-\Hat\bfB_1)\times y_\zeta\bigr)
b_1\bfB_2
- (\bfB_2\times y_\zeta) {b_1(\Hat\bfB_1-\bfA)}\\[3mm]
\qquad\qquad\qquad\qquad\qquad\qquad\quad \ds
+\bigl((\bfA-\Hat\bfB_1)\times y_\zeta\bigr) {b_2(\Hat\bfB_1-\bfA)}
-\bigl(t-\tau(\zeta)\bigr)\, (\bfB_2\times y_\zeta) \, b_2 \bfB_2\bigg\},
\enda
 \eeq
 In the above formula, only $\bfA$ contains the term $S=\sigma_t$. Hence, the fourth term in (\ref{lg0}) does not contain $S$. Also, when evaluating at $t=\tau(\zeta)$, one checks
 that the derivative of third term w.r.t. $S$ equals zero.
 Thus, we only need to compute the first two terms in (\ref{lg0}),
 namely
 $$\bega{rl} J_1&\doteq ~2 l_1^*\cdot(I+tD_x{\bfw})^{-1} \bigl((\bfA-\Hat\bfB_1)\times y_\zeta\bigr)
b_1\bfB_2,\\[2mm] J_2&\doteq\, -2 l_1^*\cdot(I+tD_x{\bfw})^{-1} (\bfB_2\times y_\zeta) {b_1(\Hat\bfB_1-\bfA)}.\enda$$
Indeed, at $t=\tau(\zeta)$,
 \begin{eqnarray*}
 \frac{\partial J_1}{\partial S}&=&2l_1^*\cdot(I+tD_x{\bfw})^{-1}\bfB_2\,b_1\,\Bigl[\Bigl(\bigl(I+tD_x\bfw\bigr)x_{\sigma}\Bigr)\times y_{\zeta}\Bigr]\\[2mm]
 &=&2l_1^*\cdot(I+tD_x{\bfw})^{-1}\sqrt{c_1(\zeta)}\lambda_1(\zeta)r_1(\zeta)b_1\,\Bigl[\Bigl(\bigl(I+tD_x\bfw\bigr)x_{\sigma}\Bigr)\times y_{\zeta}\Bigr]\\[2mm]
 &=&2\sqrt{c_1(\zeta)}\,b_1\lambda_1(\zeta)r_1\times y_{\zeta}\,.
 \end{eqnarray*}
 Here we have again used (\ref{di1})  and the fact that, at $t=\tau(\zeta)$,
 $$(I+tD_x\bfw)x_{\sigma}=(t-\tau(\zeta))\lambda_1(\zeta)r_1(\zeta).$$
In addition, at $t=\tau(\zeta)$, we have
 \begin{eqnarray*}
 \frac{\partial J_2}{\partial S}&=&\frac{\partial}{\partial S}\left[2l_1^*\cdot(I+tD_x{\bfw})^{-1}\bfA\,b_1\,\bigl(\bfB_2\times y_{\zeta}\bigr)\right]\\[3mm]
 &=&\frac{\partial}{\partial S}\left[2l_1^*\cdot\Bigl(x_{\xi}\xi_t+x_{\sigma}S\Bigr)b_1\,\Bigl(\sqrt{c_1(\zeta)}\lambda_1(\zeta)r_1\times y_{\zeta}\Bigr)\right]\\[3mm]
 &=&2\sqrt{c_1(\zeta)}\,b_1\lambda_1(\zeta)r_1\times y_{\zeta}.
 \end{eqnarray*}
 Based on the above analysis, we conclude (\ref{f1s}).

By (\ref{eig}),  the eigenvalues of ${\bf\Lambda}(\zeta)$ are found to be
 \bel{eigvv}\lambda_1~=~\left.\frac{\partial F_1}{\partial S}\right|_{t=\tau(\zeta)}~=~-3~<~-1~=~\lambda_{2,3}~<~-\frac12~=~\lambda_4~<~0~=~\lambda_5\,.\eeq

According to (\ref{eigvv}), none of the eigenvalues of $\Lambda$ is a positive integer.
Using Theorem~2.2 in \cite{bao}, or the result by Li \cite{li}, we thus
obtain  the local  existence of an analytic solution for (\ref{sip}),  defined in a neighborhood
of $t=\tau(\zeta)$ and $\zeta=0$.
This achieves the proof of Theorem~\ref{t:51}.
\endproof
%

}

\v
{\bf Acknowledgments.} The research by the first author
 was partially supported by NSF with
grant  DMS-2306926, ``Regularity and Approximation of Solutions to Conservation Laws". The work
by the second author was partially supported by NSF with grants DMS-2008504 and
DMS-2306258. The work of third author was in part supported by Zhejiang Normal University with
grants YS304222929 and ZZ323205020522016004.

\v


\begin{thebibliography}{999}

\bibitem{A} S.~Alinhac, {\it Blowup for Nonlinear Hyperbolic Equations.}
Birkh\"auser, Boston, 1995.

\bibitem{AR} A.~I.~Aptekarev and Yu.~G.~Rykov,
Emergence of a hierarchy of singularities
in zero-pressure media. Two-dimensional case.
{\it Mathematical Notes} {\bf 112} (2022), 495--504.


\bibitem{bao} M. S. Baouendi and C. Goulaouic, Singular nonlinear Cauchy problems,
{\it J.  Diff. Equat.} {\bf 22} (1976) 268-291.

%
\bibitem{Berger}
M.~S.~Berger,
{\it Nonlinearity and Functional Analysis.}  Academic Press, New York, 1977.

\bibitem{BD} S.~Bianchini and S.~Daneri, On the sticky particle solutions to the multi-dimensional pressureless Euler equations, {\it J. Differ. Equations} {\bf 368} (2023) 173-202

%
%

 \bibitem{BC} A.~Bressan and G.~Chen,
Generic regularity of conservative solutions
to a nonlinear  wave equation,
{\it Ann. Inst. Poincar\'e Anal. Nonlin.} {\bf 34} (2017), 335--354.


 \bibitem{BCH} A.~Bressan, G.~Chen, and S.~Huang, Generic singularities of 2D
 pressureless flow. {\it Science China Math.}, to appear.   Available at arXiv:2307.11602.

 \bibitem{BHY} A.~Bressan, T.~Huang, and F.~Yu,  Structurally stable singularities
for a nonlinear   wave equation. {\it Bull. Inst. Math. Acad. Sinica},
{\bf 10} (2015), 449--478.

 \bibitem{BN} A.~Bressan and T.~V.~Nguyen, Non-existence and non-uniqueness for multidimensional
 sticky particle systems.   {\it Kinetic Rel. Models} {\bf 7} (2014),
 205--218.

\bibitem{BSV}
T.~Buckmaster, S.~Shkoller, and V.~Vicol,
Shock formation and vorticity creation for 3d Euler. {\it  Comm. Pure Appl. Math.} {\bf 76} (2023), 1965-2072.

\bibitem{CO} R.~Caflisch and O.~Orellana,
Singular solutions and ill-posedness for the evolution of vortex sheets.
{\it SIAM J. Math. Anal.} {\bf  20} (1989), 29--307.

 \bibitem{CCD} H.~Cai, G.~Chen and Y.~Du,
Uniqueness and regularity of conservative solution to a wave system modeling nematic liquid crystal, {\it J. Math. Pures Appl.} {\bf  117} (2018), 185--220.

\bibitem{CH} S.-N.~Chow and J.~Hale, {\it Methods of Bifurcation Theory.} Springer, New York, 1982.
%

\bibitem{Damon} J.~Damon, Generic properties of solutions to partial differential equations
{\it Arch. Rational Mech. Anal.} {\bf 140} (1997), 353--403.

\bibitem{D} V.~G.~Danilov, Interaction of $\delta$-shock waves in a system of pressureless gas dynamics equations.
{\it Russ. J. Math. Phys.} {\bf 26} (2019),  306--319.


\bibitem{DM} V.~G.~Danilov and  D.~Mitrovic,
Shock wave formation process for a multidimensional scalar conservation law.
{\it Quart. Appl. Math.} {\bf 69} (2011), 613--634.

\bibitem{Deimling}
K.~Deimling,
{\it Nonlinear Functional Analysis.} Springer-Verlag, Berlin, 1985.
%

\bibitem{ERS} W.~E, Yu.~G.~Rykov, and  Ya.~G.~Sinai,
Generalized variational principles, global weak solutions and behavior with random initial data for systems of conservation laws arising in adhesion particle dynamics.
{\it Comm. Math. Phys.} {\bf 177} (1996), 349--380.

\bibitem{Evans}L.~C.~Evans, {\it Partial Differential Equations}. Second Edition. AMS, Providence, R.I., 2010.

%

%
%


\bibitem{GG}  M.~Golubitsky and  V.~Guillemin,
{\it Stable Mappings and Their Singularities.} Springer-Verlag, New York, 1973.

\bibitem{GS} M.~Golubitsky and D.~Schaeffer,
 Stability of shock waves for a single conservation law.
{\it Adv. Math.} {\bf 15} (1975),  65-71.


\bibitem{GSS} S.~N.~Gurbatov, A.~I.~Saichev, and S.~F.~Shandarin, Large-scale structure of the Universe. The Zeldovich
approximation and the adhesion model. {\it  Phys. Usp.}  {\bf 55} (2012),  223--249.

\bibitem{Guck}  J.~Guckenheimer, Catastrophes and Partial Differential Equations.
{\it Ann. Inst. Fourier} {\bf 23} (1973), 31--59.


\bibitem{Gk2}
J.~Guckenheimer, Solving a single conservation law. In
``Dynamical systems - Warwick 1974'', pp. 108--134.
Springer Lecture Notes in Math, {\bf 468}.

\bibitem{Kong}
D.X.~Kong,
Formation and propagation of singularities for $2\times 2$ quasilinear hyperbolic systems.
{\it Trans. Amer. Math. Soc.} {\bf  354}  (2002), 3155--3179.
%
%

\bibitem{li}F.B. ~Li, On systems of partial differential equations of
Briot-Bouquet type. {\it Bull. Aust. Math. Soc.} {\bf  98}  (2018), 122--133.

\bibitem{RR} M.~Renardy and R.~Rogers, {\it An Introduction to Partial Differential Equations.}  Second edition. Springer-Verlag, New York, 2004.

\bibitem{S} D.~Schaeffer, A regularity theorem for conservation laws,
{\it Adv. in Math.} {\bf  11}
(1973), 368--386.

\bibitem{Sever}
M.~Sever, An existence theorem in the large for zero-pressure gas dynamics.  {\it Diff. Integral Equat.}  \textbf{14}  (2001), 1077--1092.

\bibitem{ST}
R.~Shvydkoy and E.~Tadmor,  Topologically based fractional diffusion and emergent dynamics with short-range interactions. {\it SIAM J. Math. Anal.}
{\bf 52} (2020), 5792--5839.
%


\bibitem{U}  S.~Ushiki, Unstable manifold of analytic dynamical systems,
{\it J. Math Kyoto Univ.} {\bf 21} (1981), 763--785.

\bibitem{W} W.~Walter,
An elementary proof of the Cauchy-Kowalevsky theorem. {\it American Math. Monthly}
{\bf 92} (1985), 115--126.

\bibitem{Z} Ya.~B.~Zeldovich, Gravitational instability: An approximate theory for large density perturbations, {\it Astron. Astrophys.} {\bf 5} (1970), 84--89.


\end{thebibliography}
\end{document}